\input ./style/arxiv-general.cfg
\documentclass[MSNbibl,number,citesort,seceqn,dvips]{arxbj}
\makeatletter
   \@ifpackageloaded{graphicx}{}{\usepackage{graphicx}}
\makeatother
\usepackage{mathrsfs}
\usepackage{mathbh}
\usepackage{upgreek}

\volume{23}
\issue{1}
\pubyear{2017}
\firstpage{645}
\lastpage{669}
\doi{10.3150/14-BEJ677}
\docsubty{FLA}

\makeatletter
\newcommand{\ppi}{\pi}
\def\mathbbold{\mathbh}
\newcommand{\eqref}[1]{(\ref{#1})}
\newtheorem{thmm}{Theorem}[section]
\newtheorem{cor}[thmm]{Corollary}
\newtheorem{lem}[thmm]{Lemma}

\newremark{exa}[thmm]{Example}
\newproclaim{defn}{Definition}[section]
\newcommand{\scr}[1]{\mathscr #1}

\newremark{rem}{Remark}[section]

\def\R{\mathbb R}
\def\kk{\kappa}
\def\dd{\delta}
\def\vv{\varepsilon}
\def\rr{\rho}
\def\GG{\Gamma}
\def\nn{\nabla}
\def\bb{\beta}
\def\aa{\alpha}
\def\si{\sigma}
\def\F{\scr F}
\def\OO{\Omega}
\def\tt{\tilde}
\def\P{\mathbb P}
\def\E{\mathbb E}
\def\Q{\mathbb Q}
\def\LL{\Lambda}
\def\to{\rightarrow}
\def\S{\mathbb{S}}
\renewcommand{\bar}{\overline}
\renewcommand{\hat}{\widehat}
\renewcommand{\tilde}{\widetilde}

\makeatother

\begin{document}
\begin{frontmatter}

\title{Two-time-scale stochastic partial differential equations
driven by $\alpha$-stable noises: Averaging principles}

\runtitle{Two-time-scale SPDEs with $\alpha$-stable noises}

\begin{aug}
\author[A]{\inits{J.}\fnms{Jianhai}~\snm{Bao}\thanksref{A}\ead[label=e1]{jianhaibao13@gmail.com}},
\author[B]{\inits{G.}\fnms{George}~\snm{Yin}\corref{}\thanksref{B}\ead[label=e2]{gyin@math.wayne.edu}}
\and
\author[C]{\inits{C.}\fnms{Chenggui}~\snm{Yuan}\thanksref{C}\ead[label=e3]{C.Yuan@swansea.ac.uk}}
\address[A]{Department of Mathematics, Central South University,
Changsha, Hunan, 410083, P.R. China.\\ \printead{e1}}
\address[B]{Department of Mathematics, Wayne State
University, Detroit, MI 48202, USA.\\ \printead{e2}}
\address[C]{Department of
Mathematics, Swansea University, Singleton Park, SA2 8PP, UK.\\ \printead{e3}}
\end{aug}

%
\received{\smonth{7} \syear{2013}}
%
\revised{\smonth{3} \syear{2014}}

%
\begin{abstract}
This paper focuses on stochastic partial differential
equations (SPDEs) under two-time-scale formulation. Distinct from
the work in the existing literature, the systems are
driven by $\aa$-stable processes with $\aa\in(1,2)$.
In addition, the SPDEs are either modulated by a continuous-time
Markov chain with a finite state space or have an addition fast jump
component. The inclusion of the Markov chain is for the
needs of treating random environment, whereas the addition of the
fast jump process enables the consideration of
discontinuity in the sample paths of the fast processes.
Assuming either a fast changing Markov switching or an additional
fast-varying jump process,
this work aims to
obtain the averaging principles for such systems. There are several
distinct difficulties. First, the noise is not square integrable.
Second, in our setup, for the underlying SPDE, there is only a
unique mild solution and as a result, there is only mild It\^o's
formula that can be used. Moreover, another new aspect is the
addition of the fast regime switching and the addition of the fast
varying jump processes in the formulation, which enlarges the
applicability of the underlying systems. To overcome these
difficulties, a semigroup approach is taken. Under suitable
conditions, it is proved that the $p$th moment convergence takes
place with $p \in(1,\aa)$, which is stronger than the usual weak
convergence approaches.
\end{abstract}

%
\begin{keyword}
\kwd{$\aa$-stable process}
\kwd{averaging principle}
\kwd{invariant measure}
\kwd{stochastic partial differential equation}
\kwd{strong convergence}
\end{keyword}
\end{frontmatter}

\section{Introduction}
Averaging principles for stochastic differential equations (SDEs)
have been studied extensively, for example, in
Liu and Vanden-Eijnden \cite{ELV},
Freidlin and Wentzell \cite{FW},
Khasminskii \cite{K68},
Yin and Zhang \cite{YZ}.
Recently, averaging principles for stochastic partial differential
equations (SPDEs) have also drawn
much attention; see, for example,
Kuksin and Piatnitski \cite{KP} and
Maslowski \textit{et al.} \cite{MSV}. In particular,
Bl\"omker \textit{et al.} \cite{BHP} derived averaging results with explicit
error bounds for
SPDEs with quadratic nonlinearities, where the limiting system is an
SDE;
Cerrai and Freidlin \cite{CF} investigated the weak convergence
for two-time-scale stochastic reaction--diffusion equations with
additive noise by using an approach based on Kolmogorov equations
and martingale solutions of stochastic equations;
Cerrai \cite{C09}
generalized Cerrai and Freidlin \cite{CF} to the case of slow--fast
reaction--diffusion
equations driven by multiplicative noise, where
the reaction terms appear in both equations;
Br\'{e}hier \cite{B12} gave
the strong and weak orders in averaging for stochastic evolution
equation of parabolic type with slow and fast time scales. For the
finite-dimensional jump--diffusion case, we refer to
Givon \cite{G07}.

In view of the development on the aforementioned singularly
perturbed SPDEs, the noise processes considered to date are
mainly
square integrable processes. However,
such
requirement rules out the interesting $\aa$-stable processes.
It is well known that both Wiener processes and Poisson-jump processes have
finite moments of any order,
whereas an $\aa$-stable process only has finite $p$th
moment for $p\in(0,\aa)$.
Stochastic equations driven by
$\aa$-stable processes
have proven to have numerous applications in physics
because such processes can be used to model systems with
heavy
tails. As a result, such processes
have received
increasing attentions recently. For example,
Priola and Zabczyk \cite{PZ11} gave a proper starting point on the
investigation of structural properties of SPDEs driven by an
additive cylindrical stable noise;
Dong \textit{et al.} \cite{DXZ} studied
ergodicity of stochastic Burgers equations driven by
$\aa/2$-subordinated cylindrical Brownian motions with
$\aa\in(1,2)$.
For
finite-dimensional SDEs driven by $\aa$-stable noises,
Wang \cite{W11a} derived gradient estimate for linear SDEs,
Zhang \cite{Z12} established the Bismut--Elworthy--Li derivative formula for
nonlinear SDEs, and
Ouyang \cite{O09}
established Harnack
inequalities for Ornstein--Uhlenbeck processes by the sharp estimates
of density function for rotationally invariant symmetric
$\aa$-stable L\'{e}vy processes. Nevertheless, two-time-scale
formulation for stochastic processes driven by $\aa$-stable
processes have not yet been considered to date to the best of our
knowledge.

Motivated by the previous works, in this paper we
develop
averaging principles for two-time-scale SPDEs driven by $\aa$-stable
noises
that admit unique mild solutions. The time-scale separation is given
by introducing a small parameter $\vv>0$. For the case of
mean-square integrable noise, the It\^o formula plays an important
role in the error analysis between the slow component and the
averaging systems;
see, for example,
Givon \cite{G07},
Fu and Duan \cite{FD}
and
Fu and Liu \cite{FL}.
It has been noted that when the diffusion operators in Fu and Duan
\cite{FD} and
Fu and Liu \cite{FL} are Hilbert--Schmidt, the mild solution is indeed
a strong
solution. Nevertheless,
in our case, only \textit{mild It\^o's formula} (see, e.g.,
Da Prato \textit{et al.} \cite{DJR}, Theorem~1) is available since the
stochastic systems
considered only admit mild solutions, not strong solutions.
Moreover, the technique adopted in Br\'{e}hier \cite{B12}, Lemma~3.1,
which is a
key ingredient in discussing averaging principle, does not work for
the case of SPDEs driven by $\aa$-stable noises either, although the
mild solution is treated there. In our study, in addition to the
SPDEs, we assume that the systems are modulated by a continuous-time
Markov chain. This Markov chain has a finite state space resulting
in a system of stochastic differential equations switching back and
forth according to the state of the Markov chain. The Markov chain
can be used to model discrete events that are not representable
otherwise. It is by now widely recognized that such regime-switching
formulation is an effective way of modeling many practical
situations in which random environment and other random factors have
to be taken into consideration. Perhaps, one of the first efforts in
modeling random environment using a finite-state Markov chain can be
traced back to
Griego and Hersh \cite{GriegoH69} (see also the
extended survey in Hersh \cite{Hersh74}, where multiple time scale was
also used). Much of the recent modeling and analysis effort stems
from the work of
Hamilton and Susmel \cite{HamS}, who revealed the feature of the
so-called regime-switching
systems under which the dynamics of the systems can be quite
different under different regimes. Their idea stimulated much of the
subsequent study. For example, in the simplest setting, the
successfully used regime-switching models in financial market
portraits the random environment with two states bull and bear
markets, whose volatilities are drastically different.

Our study is divided into two parts.
In the first part, we assume that the switching process is subject
to fast variation, either within a weakly irreducible class or
within a number of nearly decomposable weakly irreducible classes
(see Yin and Zhang \cite{YZ}, Chapter~4). The idea is that the
original system
subject to fast switching is more complex, but the limit system is
much simpler. For many applications, it will be desirable to find
the structure of the limit system leading substantial reduction of
computational complexity for such tasks as control and optimization
etc. We show that under suitable conditions, a limit process that is
a solution of either an SPDE or an SPDE with switching is obtained.
The key is that in the limit, the coefficients are averaged out with
respect to the stationary measure of the switching processes. In the
second part,
we assume that there is an
additional fast-varying
random process. Although the process is fast varying, it does not
blow up, but rather has an invariant measure. The ergodicity of the
fast process helps us to get a limit process that is a solution of
the SPDEs with the coefficients being averaged out with respect to
the stationary distribution of the fast-varying process.

To summarize, there are several distinct difficulties in our
problems. First, the noise is not square integrable. Second, the
underlying SPDE admits only a unique mild solution and as a result,
there is only mild It\^o's formula that can be used. Moreover,
another new aspect is the addition of the fast regime switching and
the addition of the fast varying jump processes in the formulation,
which enlarges the applicability of the underlying systems. To
overcome these difficulties, using the mild solutions, a semigroup
approach is taken. Under suitable conditions, it is proved that the
$p$th moment convergence takes place with $p \in(1,\aa)$, which is
stronger than the usual weak convergence approaches. We thus term
such a convergence as strong convergence.

The rest of the paper is organized as follows. In Section~\ref{sec:2-sw},
we obtain not only averaging principles for SPDEs
with two-time-scale Markov switching with a single weakly recurrent
class but also for the case of two-time-scale Markov switching with
multiple weakly irreducible classes. In Section~\ref{fast-jump}, we
demonstrate the strong convergence for SPDEs with an additional
fast-varying random process driven by cylindrical stable processes.

\section{SPDEs with two-time-scale Markov switching}\label{sec:2-sw}
We first
recall some basics on
stable processes. A real-valued random
variable $\eta$ is said to have a stable distribution with
stability index $\aa\in(0,2)$, scale parameter $\si\in(0,\infty )$,
skewness parameter $\bb\in[-1,1]$, and location parameter
$\mu\in(-\infty,\infty )$, if
its characteristic function has the
form:
\[
\phi_{\eta}(u)=\E\exp(iu\eta)=\exp\bigl\{-\si^\aa|u|^\aa\bigl(1-i
\beta \operatorname {sgn}(u)\Phi\bigr)+i\mu u\bigr\},\qquad  u\in\R,
\]
where $\Phi=\tan(\uppi\aa/2)$ for $\aa\neq1$ and
$\Phi=-(2/\uppi)\log|u|$ for $\aa=1$.
Note that
the monograph Samorodnitsky and Taqqu \cite{ST94}, pages~2--10, also gives
three other
equivalent definitions of a stable distribution. We denote the
family of stable distributions by $S_\aa(\si,\bb,\mu)$ and write
$X\sim S_\aa(\si,\bb,\mu)$ to indicate that $X$ has the stable
distribution $S_\aa(\si,\bb,\mu)$. A random variable $X\sim
S_\aa(\si,\bb,\mu)$ is said to be strictly stable if $\mu=0$ for
$\aa\neq1$ (if $\bb=0$ for $\aa=1$),
symmetric if $\bb=\mu=0$,
and
standard (normalized)\vadjust{\goodbreak} if $\bb=\mu=0$ and $\si=1$.
Let
$(H,\langle \cdot,\cdot\rangle ,|\cdot|_H)$ be a real separable
Hilbert space.
Let
$L=\{L(t)\}_{t\ge0}$ and $Z=\{Z(t)\}_{t\ge0}$ be a cylindrical
$\aa$-stable process and $\bb$-stable process defined by the
orthogonal expansion, respectively,
%
\begin{equation}
\label{a1} L(t):=\sum_{k=1}^\infty
\bb_kL_k(t)e_k \quad\mbox{and}\quad Z(t):=\sum
_{k=1}^\infty q_kZ_k(t)e_k,
\end{equation}
where $\{e_k\}_{k\ge1}$ is an orthonormal basis of $H$,
$\{L_k(t)\}_{k\ge1}$ and $\{Z_k(t)\}_{k\ge1}$ are sequences of
i.i.d. (independent and identically distributed) real-valued
symmetric $\aa$-stable processes and $\bb$-stable processes defined
on the stochastic basis $(\OO, \F, {\{\F_t}\}_{t\ge0}, \P)$,
respectively, and $\bb_k,q_k>0$ for each $k\ge1$. $\|\cdot\|$
stands for the operator norm, and $\mathscr{L}(\xi)$ means the law
of an $H$-valued random variable $\xi$. Throughout this paper, we
assume that $\aa,\bb\in(1,2]$. Generic constants will be denoted by
$c$, and we use the shorthand notation $a\lesssim b$ to mean $a\le
cb$. If the constant $c$ depends on a parameter $p$, we shall also
write $c_p$ and $a\lesssim_pb$.

\subsection{
Two-time-scale Markov switching with a single weakly irreducible
class}

Hybrid systems driven by continuous-time Markov chains have been
used to model many practical scenarios in which abrupt changes may
be experienced in the structure and parameters caused by phenomena
such as component failures or repairs; see Sethi and Zhang \cite{SethiZ},
Remark~5.1, page~94,
for discussions on the modeling of such a system
and related optimal control problems. For finite-dimensional cases,
there is extensive literature on such topic, for example,
Mao and Yuan \cite{MY06},
Mariton \cite{M90},
Yin and Zhu \cite{YZ10}
and the
references therein. As an infinite-dimensional
example,
we consider a one-dimensional rod of length $\uppi$
whose ends are maintained at $0^\circ$ and whose sides are
insulated. Assume that there is an exothermic reaction taking place
inside the rod with heat being produced proportionally to the
temperature. The temperature $u$ in the rod may be modeled by
%
\begin{equation}
\label{eq111} %
\cases{ \displaystyle\frac{\partial u}{\partial t}=
\frac{\partial^2u}{\partial x^2}+cu,&\quad $t>0, x
\in(0,\uppi),$\vspace*{2pt}
\cr
u(t,0)=u(t,\uppi)=0,&\quad $u(0,x)=u_0(x),$}
\end{equation}
where $u=u(t,x)$ and $c$ is a constant dependent on the rate of
reaction.
In lieu of assuming the system to be in a fixed configuration, let
system \eqref{eq111}
switch from one mode to
another in a random way when it experiences abrupt changes in its
structure and parameters caused by phenomena such as component
failures or repairs, changing subsystem interconnections, or abrupt
environmental disturbances.
The system under regime switching could be described by the
following random model
\[
\cases{ \displaystyle\frac{\partial u}{\partial t}=\frac{\partial^2u}{\partial
x^2}+c\bigl(r(t)\bigr)u,&\quad $t>0, x
\in(0,\uppi),$\vspace*{2pt}
\cr
u(t,0)=u(t,\uppi)=0,&\quad $u(0,x)=u_0(x),
r(0)=r_0,$} %
\]
where $r(t)$ is a
continuous-time Markov chain with a finite state space $\mathbb{S}$
and $c\dvtx \mathbb{S}\rightarrow\mathbb{R}$. For further details, we
refer to, for example,
Anabtawi \cite{A11} and
Bao \textit{et al.} \cite{BMY}.

With the motivation above, assuming that $r^\vv(t)$ is a
continuous-time Markov chain with a finite state space
$\S:=\{1,\ldots,n\}$,
we
consider the following SPDE
%
\begin{equation}
\label{eq1} \mathrm{d}X^\vv(t)=\bigl\{AX^\vv(t)+b
\bigl(X^\vv(t),r^\vv(t)\bigr)\bigr\}\,\mathrm{d}t+
\mathrm{d}L(t),\qquad t>0
\end{equation}
with the initial values $X^\vv(0)=x\in H$ and $r^\vv(0)=r_0\in\S$.

In \eqref{eq1}, for any $\vv\in(0,1)$, $r^\vv(t)$ is a Markov
chain with a finite state space $\S$ and generator
\[
\Q^\vv:={\tt\Q\over\vv}+\hat\Q,
\]
where $\tt\Q$ and $\hat\Q$ are
suitable generators of some
Markov chains such that
$\tt\Q/\vv$ and $\hat\Q$ represent the
fast-varying and the slow-changing parts, respectively.
In what follows, we further assume that $\tt\Q$ is weakly
irreducible.
That is, the system of equations
%
\begin{equation}
\label{eq29} \cases{
\nu\tt\Q=0,
\vspace*{2pt}\cr
\displaystyle\sum_{i=1}^n
\nu_i=1, }
\end{equation}
has a unique solution satisfying $\nu_i \ge0$ for all $i\in{\mathbb{S}}$.
Throughout this subsection, we
assume that the following conditions fulfill.
\begin{longlist}[(A1)]
\item[(A1)] $A\dvtx \mathcal{D}(A)\subset H\mapsto H$ is a
self-adjoint compact operator on $H$ such that $-A$ has
discrete spectrum
$0<\lambda_1<\lambda_2<\cdots<\lambda_k<\cdots$ and $\lim_{k\to\infty }\lambda_k=\infty $. In
this case, $A$ generates an
analytic contraction
semigroup $\{\mathrm{e}^{tA}\}_{t\ge0}$, such that $\|\mathrm{e}^{tA}\|\leq
\mathrm{e}^{-\lambda_1t}$.
\item[(A2)] For each $i\in\S$, there exists $ K_i>0$ such that
\[
\bigl|b(x,i)-b(y,i)\bigr|_H\le K_i|x-y|_H,\qquad x,y
\in H.
\]
\item[(A3)] There exists $\theta\in(0,1)$ such that
$\aa\theta\in(0,1)$ and
\[
\dd:=\sum_{k=1}^\infty \frac {\bb_k^\aa}{\lambda_k^{1-\aa\theta
}}<
\infty.
\]
\end{longlist}

Under assumption (A1)--(A3), according to Mao and Yuan \cite{MY06},
Theorem~3.13, page~89, and
Priola and Zabczyk \cite{PZ11}, Theorem~5.3, \eqref{eq1} admits a
unique mild solution,
that is, there exists a predictable $H$-valued stochastic process
$\{X^\vv(t)\}_{t\ge0}$ such that
%
\begin{equation}
\label{mild-sl} X^\vv(t)= \mathrm{e}^{ A t} x + \int
^t_0 \mathrm{e}^{
A(t-s)} b
\bigl(X^\vv(s),r^\vv(s)\bigr)\,\mathrm{d} s+\int
^t_0 \mathrm{e}^{ A(t-s)}\,\mathrm{d}L(s),\qquad
\P\mbox{-a.s.}
\end{equation}
As can be seen,
compared with the fast varying
$r^\vv(\cdot)$,
$X^\vv(\cdot)$ changes relatively slowly.
The intuitive idea can be explained as follows.
Using
the methods of stochastic averaging initiated in Khasminskii \cite
{K66} (see also
Khasminskii \cite{K68}, Khasminskii and Yin \cite{KhY04,KhY05})
and subsequently developed by
Kushner \cite{Kushner84},
$r^\vv(t)$ can be treated essentially as a ``noise'' process.
With the slow variable ``fixed'' or ``frozen,''
a law of large numbers holds so the noise is averaged out.
Moreover,
the slow component $X^\vv(t)$ converges to $\bar X(t)$ in an
appropriate
sense.
We will show that
the limit $\{\bar X(t)\}_{t\ge0}$ satisfies in the mild sense an
SPDE
%
\begin{equation}
\label{eq8} \mathrm{d}\bar X(t)=\bigl\{A\bar X(t)+\bar b\bigl(\bar X(t)\bigr)
\bigr\}\,\mathrm{d}t+\mathrm{d}L(t),\qquad t>0,
\end{equation}
with initial value $\bar X(0)=x\in H$, where $\bar
b(x):=\sum_{i=1}^nb(x,i)\nu_i,x\in H$, an average with respect to
the invariant measure $\nu:=(\nu_1,\ldots,\nu_n)$ given in
\eqref{eq29}. Our main
result of this section is
given as follows.

\begin{thmm}\label{two}
Let \textup{(A1)--(A3)} hold and assume further that the initial value
$X^\vv(0)=x\in\mathcal{D}((-A)^\theta)$.
Then, for any sufficiently small $\vv\in(0,1)$,
\[
\sup_{0\le t\le T}\bigl(\E\bigl|X^\vv(t)-\bar
X(t)\bigr|_H^p\bigr)^{1/p}\lesssim _T
\vv^{\rr\theta},\qquad p\in(1,\aa),
\]
where $\theta\in(0,1)$ is the constant such that \textup{(A3)} holds and
$\rr<(\aa-p)/(\aa-p+p\theta\aa)$.
\end{thmm}

To prove
Theorem~\ref{two}, we need the following lemma.

\begin{lem}\label{Le}
Let the assumptions of Theorem~\ref{two} hold. Then, for any
$h\in(0,1)$ and $p\in(1,\aa)$,
\[
\sup_{0\le t\le T}\bigl(\E\bigl|\bar X({t+h})-\bar X(t)\bigr|_H^p
\bigr)^{1/p}\lesssim _Th^{\theta}. 
\]
\end{lem}

\begin{pf}
Noting that $(\E|\cdot|_H^p)^{1/p}, p\in(1,\aa)$, is a norm,
we get from (A1), (A2), and Priola and Zabczyk \cite{PZ11}, (4.12), that
%
\begin{eqnarray}
\label{A1} %
&&\bigl(\E\bigl|\bar X(t)\bigr|_H^p
\bigr)^{1/p}\nonumber
\\
&&\quad\le|x|_H+\int_0^t\bigl\|
\mathrm{e}^{(t-s)A}\bigr\|\bigl(\E\bigl|\bar b\bigl(\bar X(s)\bigr)\bigr|_H^p
\bigr)^{1/p}\,\mathrm{d}s+\E \biggl(\biggl |\int_0^t
\mathrm{e}^{(t-s)A}\,\mathrm{d} L(s) \biggr|_H^p
\biggr)^{1/p}
\nonumber
\\[-8pt]
\\[-8pt]
\nonumber
&&\quad\le|x|_H+\sum_{i=1}^n
\nu_i\int_0^t
\mathrm{e}^{-\lambda
_1(t-s)}\bigl(\E\bigl| b\bigl(\bar X(s),i\bigr)\bigr|_H^p
\bigr)^{1/p}\,\mathrm{d} s+c \Biggl(\sum_{k=1}^\infty
\frac {\bb_k^\aa(1-\mathrm{e}^{-\aa
\lambda_kt})}{\lambda_k} \Biggr)^{1/\aa}\quad
\\
&&\quad\le|x|_H+\sum_{i=1}^n
\nu_i\int_0^t
\mathrm{e}^{-\lambda
_1(t-s)}\bigl\{K_i\bigl(\E\bigl| \bar
X(s)\bigr|_H^p\bigr)^{1/p}+\bigl(\E\bigl|
b(0,i)\bigr|_H^p\bigr)^{1/p}\bigr\}\,\mathrm{d} s+c
\tau,
\nonumber
\end{eqnarray}
where $\tau:= (\sum\frac {\bb_k^\aa}{\lambda_k} )^{1/\aa
}<\infty $
according to (A3). Multiplying $\mathrm{e}^{\lambda_1t}$ on both
sides of
\eqref{A1} gives
\[
\mathrm{e}^{\lambda_1t}\bigl(\E\bigl|\bar X(t)\bigr|_H^p
\bigr)^{1/p}\le c\bigl(1+\mathrm{e}^{\lambda_1
t}\bigr)+\sum
_{i=1}^n\nu_iK_i\int
_0^t\mathrm{e}^{\lambda_1s}\bigl(\E\bigl| \bar
X(s)\bigr|_H^p\bigr)^{1/p}\,\mathrm{d}s.
\]
This, together with the Gronwall inequality, yields that
%
\begin{equation}
\label{eq15} \sup_{0\le t\le T}\E\bigl|\bar X(t)\bigr|_H^p<
\infty.
\end{equation}
From \eqref{eq8}, one has
%
\begin{eqnarray}
\label{eq18} %
&&\bigl(\E\bigl|\bar X(t+h)-\bar X(t)\bigr|_H^p
\bigr)^{1/p}\nonumber
\\
&&\quad\le\bigl|\bigl(\mathrm{e}^{hA}-{\mathbf1}\bigr)\mathrm{e}^{t
A}x\bigr|_H+
\sum_{i=1}^n\nu_i\int
_0^t\bigl(\E\bigl|\bigl(\mathrm{e}^{hA}-{
\mathbf 1}\bigr)\mathrm{e}^{(t-s)A}b\bigl(\bar X(s),i\bigr)
\bigr|_H^p\bigr)^{1/p}\,\mathrm{d}s
\nonumber
\\
&&\qquad{}+\sum_{i=1}^n\nu_i\int
_{ t}^{t+h}\bigl(\E\bigl|\mathrm{e}^{(t+h-s)A}b\bigl(
\bar X(s),i\bigr)\bigr|_H^p\bigr)^{1/p}\,
\mathrm{d}s
\nonumber
\\[-8pt]
\\[-8pt]
\nonumber
&&\qquad{}+ \biggl(\E \biggl|\int_0^t\bigl(
\mathrm{e}^{hA}-{\mathbf1}\bigr)\mathrm{e}^{(
t-s)A}\,
\mathrm{d}L(s)\biggr|_H^p\biggr)^{1/p}
\\
&&\qquad{} + \biggl(\E \biggl|\int_{
t}^{t+h}
\mathrm{e}^{(t+h-s)A}\,\mathrm{d}L(s) \biggr|_H^p
\biggr)^{1/p}
\nonumber
\\
&&\quad=:\sum_{j=1}^5\LL_j(t),
\nonumber
\end{eqnarray}
where ${\mathbf1}$ is the
identity operator on $H$. By the spectral properties of operator
$A$, observe that
%
\begin{equation}
\label{eq27} 
\bigl\|(-A)^\dd\mathrm{e}^{tA}\bigr\| \le
\mathrm{e}^{-\dd}\dd^\dd t^{-\dd},\qquad t>0, \dd\in(0,1)
\end{equation}
and
that
%
\begin{equation}
\label{eq28} %
\bigl\|(-A)^{-\delta}\bigl({\mathbf1}-
\mathrm{e}^{tA}\bigr)\bigr\| \le c_\dd t^\dd,\qquad  t>0,
\dd\in(0,1)
\end{equation}
for some $c_\dd>0$.
Due to $x\in\mathcal{D}((-A)^\theta)$, taking (A1), (A2),
\eqref{eq15}, \eqref{eq27}, and \eqref{eq28} into account yields
that
%
\begin{eqnarray}
\label{eq20} %
\LL_1(t)+\LL_2(t)&\le&\bigl\|\bigl(
\mathrm{e}^{hA}-{\mathbf1}\bigr) (-A)^{-\theta
}\bigr\|\cdot\bigl\|
\mathrm{e}^{t
A}\bigr\|\cdot\bigl|(-A)^{\theta}x\bigr|_H
\nonumber\\
&&{} +\sum_{i=1}^n\nu_i\int
_0^t\bigl\|\bigl(\mathrm{e}^{hA}-{\mathbf
1}\bigr) (-A)^{-\theta}\bigr\|\cdot \bigl\|\mathrm{e}^{ {(t-s)A}/{2}}\bigr\|\nonumber
\\
& &\hspace*{50pt}{}\times\bigl\|\mathrm{e}^{{(t-s)A}/{2}}(-A)^\theta\bigr\|\bigl(\E\bigl| b\bigl(\bar
X(s),i\bigr)\bigr|_H^p\bigr)^{1/p}\,\mathrm{d}s
\\
&\lesssim _T &\biggl(1+\int_0^{t}
\mathrm{e}^{- {\lambda
_1(t-s)}/{2}} \biggl(\frac {
t-s}{2} \biggr)^{-\theta}\,
\mathrm{d} s \biggr)h^\theta\nonumber
\\
&\lesssim _T&\bigl(1+\Gamma(1-\theta)\bigr)h^\theta,\nonumber
\end{eqnarray}
where $\GG(\cdot)$ is the Gamma function. Also, by (A1), (A2), and
\eqref{eq15}, we arrive at
%
\begin{equation}
\label{eq21} \LL_3(t)\lesssim _Th.
\end{equation}
Note that
\[
\int_0^{t}(-A)^{\theta}
\mathrm{e}^{( t-s)A}\,\mathrm{d} L(s)=\sum_{k=1}^\infty
\biggl(\bb_k\lambda_k^\theta\int
_0^t\mathrm{e}^{-(t-s)\lambda_k}\,\mathrm{d}
L_k(s) \biggr)e_k.
\]
Upon using an argument similar to that of Priola and Zabczyk \cite{PZ11},
Theorem~4.5,
we obtain from~(A3) that
%
\begin{eqnarray}
\label{eq12} \biggl(\E \biggl|\int_0^{t}(-A)^{\theta}
\mathrm{e}^{(
t-s)A}\,\mathrm{d} L(s) \biggr|_H^p
\biggr)^{1/p} \lesssim \Biggl(\sum_{k=1}^\infty
\bb_k^\aa\frac {1}{\aa\lambda
_k^{1-\aa\theta
}}\bigl(1-e^{-\aa
\lambda_k
t}
\bigr) \Biggr)^{1/\aa}\lesssim \dd^{1/\aa},
\end{eqnarray}
and that
%
\begin{equation}
\label{eq16} \LL_5(t)\lesssim \Biggl(\sum
_{k=1}^\infty \frac {\bb_k^\aa
}{\lambda_k}\bigl(1-
\mathrm{e}^{-\lambda_kh}\bigr) \Biggr)^{1/\aa}\lesssim \Biggl(\sum
_{k=1}^\infty \frac {\bb_k^\aa
}{\lambda_k}(\lambda
_kh)^{\aa\theta
} \Biggr)^{1/\aa}\lesssim
\dd^{1/\aa}h^\theta,
\end{equation}
where we have used the fundamental inequality: for any
$\gamma\in(0,1]$, there exists $c_\gamma>0$ such that
\[
\bigl|\mathrm{e}^{-x}-\mathrm{e}^{-y}\bigr|\le c_\gamma|x-y|^\gamma,\qquad
x,y\ge0.
\]
Thus we deduce from \eqref{eq28}, \eqref{eq12}, and \eqref{eq16}
that
%
\begin{equation}
\label{eq22} \LL_4(t)+\LL_5(t) \lesssim
_Th^\theta.
\end{equation}
As a result, the desired assertion follows by substituting
\eqref{eq20}, \eqref{eq21}, and \eqref{eq22} into \eqref{eq18}.
\end{pf}

With the aid of Lemma~\ref{Le},
we complete the proof
Theorem~\ref{two} in what follows.

\begin{pf*}{Proof of Theorem~\ref{two}} It follows from
\eqref{eq1} and \eqref{eq8} that
\begin{eqnarray*}
\bigl(\E\bigl|X^\vv(t)-\bar X(t)\bigr|_H^p
\bigr)^{1/p} &\le&\sum_{i=1}^n\int
_0^t\bigl(\E\bigl|\mathrm{e}^{(t-s)A}\bigl\{b
\bigl(X^\vv (s),i\bigr)-b\bigl(\bar X(s),i\bigr)\bigr\}
\bigr|_H^p\bigr)^{1/p}\,\mathrm{d}s
\\
&&{} +\sum_{i=1}^n \biggl(\E \biggl|\int
_0^t\mathrm{e}^{(t-s)A}b\bigl(\bar X(s),i
\bigr)\{{\mathbbold{1}}_{\{r^\vv(s)=i\}}-\nu_i\}\,\mathrm{d}s
\biggr|_H^p \biggr)^{1/p}
\\
&=:&\Xi_1(t)+\sum_{i=1}^n
\Xi_{2i}(t),
\end{eqnarray*}
where $\mathbbold{1}_\Gamma$ is the indicator function of a set $\GG$.
Taking (A1) and (A2)
into account, we have
\[
\Xi_1(t) \le \sum_{i=1}^nK_i
\int_0^t\mathrm{e}^{-\lambda_1(t-s)}\bigl(
\E\bigl|X^\vv (s)-\bar X(s)\bigr|_H^p
\bigr)^{1/p}\,\mathrm{d}s.
\]
Next, note that from the boundedness of
$|{\mathbbold{1}}_{\{r^\vv(s)=i\}}-\nu_i|$,
\begin{eqnarray*}
\Xi_{2i}(t) &\le&\int_0^t\bigl\|
\mathrm{e}^{(t-s)A}\bigl({\mathbf1}-\mathrm{e}^{(s-\lfloor
s\rfloor )A}\bigr)\bigr\|\bigl(
\E\bigl|b\bigl(\bar X(s),i\bigr)\bigr|_H^p\bigr)^{1/p}\,
\mathrm{d} s
\\
& &{}+\int_0^t\bigl\|\mathrm{e}^{(t-\lfloor s\rfloor )A}\bigr\|
\bigl(\E\bigl|b\bigl(\bar X(s),i\bigr)-b\bigl(\bar X\bigl(\lfloor s\rfloor \bigr),i
\bigr)\bigr|_H^p\bigr)^{1/p}\,\mathrm{d}s
\\
&&{} + \biggl(\E\biggl |\int_0^t\mathrm{e}^{(t-\lfloor s\rfloor
)A}b
\bigl(\bar X\bigl(\lfloor s\rfloor \bigr),i\bigr)
\{{\mathbbold{1}}_{\{r^\vv(s)=i\}}-
\nu_i\}\,\mathrm{d}s\biggr |_H^p
\biggr)^{1/p}
\\
&=:&\Upsilon_{1i}(t)+\Upsilon_{2i}(t)+\Upsilon_{3i}(t),
\end{eqnarray*}
where $\lfloor s\rfloor:=[s/\vv^\rho]\vv^\rho$ with $[s/\vv
^\rho]$
denoting the integer part of $s/\vv^\rho$ for
$\rr<(\aa-p)/(\aa-p+p\theta\aa)$. By a similar calculation as in
\eqref{eq20}, one has
%
\begin{equation}
\label{i3} %
\Upsilon_{1i}(t) \lesssim _T
\vv^{\rho\theta}. %
\end{equation}
By virtue of (A1), (A2), and Lemma~\ref{Le}, it follows that
%
\begin{equation}
\label{i4} %
\Upsilon_{2i}(t)\lesssim \int
_0^t\mathrm{e}^{-\lambda
_1(t-\lfloor s\rfloor
)}\bigl(\E\bigl|\bar
X(s)-\bar X\bigl(\lfloor s\rfloor \bigr)\bigr|_H^p
\bigr)^{1/p}\,\mathrm{d}s \lesssim _T\vv^{\rho\theta}.
\end{equation}
Let $t_j:=j\vv^\rho,j=0,\ldots,[ t/\vv^\rho]$, and $t_{[
t/\vv^p]+1}:=t$. Then, an application of the H\"older inequality
gives that
%
\begin{eqnarray}
\label{eq25} %
\Upsilon_{3i}(t) &\le&\sum
_{j=0}^{\lfloor t/\vv^\rho\rfloor } \biggl\{\E\bigl|\mathrm{e}^{(t-t_j)A}b
\bigl(\bar X(t_j),i\bigr)\bigr|_H^p \biggl|\int
_{t_j}^{t_{j+1}}\{{\mathbf 1}_{\{r^\vv(s)=i\}}-
\nu_i\}\,\mathrm{d}s \biggr|^p \biggr\}^{1/p}\nonumber
\\
&\le&\sum_{j=0}^{\lfloor
t/\vv^\rho\rfloor } \bigl(\E\bigl|
\mathrm{e}^{(t-t_j)A}b\bigl(\bar X(t_j),i\bigr)\bigr|_H^{p(1+\dd)}
\bigr)^{1/(p(1+\dd))}
\\
&&\hspace*{19pt}{} \times \biggl(\E \biggl|\int_{t_j}^{t_{j+1}}\{{\mathbf
1}_{\{r^\vv(s)=i\}}-\nu_i\}\,\mathrm{d} s \biggr|^{(p(1+\dd))/\dd}
\biggr)^{\dd/(p(1+\dd))}\nonumber %
\end{eqnarray}
for arbitrary $0<\dd<(\aa-p)/p$. Thanks to $\aa\in(1,2)$ and
$p\in(1,\aa)$, one has $p>(2\aa)/(2+\aa)$, which further gives that
$(\aa-p)/p<p/(2-p)$. Then, for $0<\dd<(\aa-p)/p$, we have
%
\begin{equation}
\label{eq14} p(1+\dd)<\aa\quad\mbox{and}\quad \bigl(p(1+\dd)\bigr)/\dd>2.
\end{equation}
Hence, (A1), \eqref{eq15}, and \eqref{eq14} yield that
%
\begin{equation}
\label{eq23} \bigl(\E\bigl|\mathrm{e}^{(t-t_j)A}b\bigl(\bar X(t_j),i
\bigr)\bigr|_H^{p(1+\dd)} \bigr)^{1/(p(1+\dd))}\lesssim
_T\mathrm{e}^{-\lambda_1(t-t_j)}.
\end{equation}
We
claim
that
%
\begin{equation}
\label{eq24} \biggl(\E \biggl|\int_{t_j}^{t_{j+1}}\{{
\mathbbold{1}}_{\{r^\vv(s)=i\}
}-\nu _i\}\,\mathrm{d} s
\biggr|^{(p(1+\dd))/\dd} \biggr)^{\dd/(p(1+\dd))}\lesssim \vv^{\rr
+((\bb-
\rr)\dd)/(p(1+\dd))},
\end{equation}
for sufficiently small $\vv\in(0,1)$, where $\bb\in(\rho,1)$ is some
constant. To show \eqref{eq24}, we adopt an argument similar to that
of Yin and Zhang \cite{YZ}, Theorem~7.2, page~170. Let
\[
\eta^\vv(u):= {1\over
2}\E \biggl|\int_{t_j}^u
\{{\mathbbold{1}}_{\{r^\vv(s)=i\}}-\nu_i\} \,\mathrm{d} s
\biggr|^2,\qquad u\in[t_j,t_{j+1}].
\]
Then, it is easy to see from the chain rule that
\[
\frac {\mathrm{d}\eta^\vv(u)}{\mathrm{d}
u}=\E\int_{t_j}^u\bigl\{({
\mathbbold{1}}_{\{r^\vv(u)=i\}}-\nu _i) ({\mathbbold
{1}}_{\{r^\vv(s)=i\}}-\nu_i)\bigr\}\,\mathrm{d} s,\qquad u
\in[t_j,t_{j+1}].
\]
Let $\bar t_k:=k\vv^\bb,k=0,1,\ldots,[(u-t_j)/\vv^\bb]+1$, where
$\bar t_0:=t_j$ and $\bar t_{\lfloor (u-t_j)/\vv^\bb\rfloor
+1}:=u$. Thus, by
the boundedness of $|{\mathbbold{1}}_{\{r^\vv(s)=i\}}-\nu_i|$,
we obtain that
\begin{eqnarray*}
\frac {\mathrm{d}\eta^\vv(u)}{\mathrm{d}u}&=&\E\int_{\bar
t_0}^{\tt
t_j}\bigl\{({
\mathbbold{1}}_{\{r^\vv(u)=i\}}-\nu_i) ({\mathbbold{1}}_{\{
r^\vv
(s)=i\}}-
\nu_i)\bigr\}\,\mathrm{d}s
\\
&&{} +\E\int_{\tt t_j}^t\bigl\{({\mathbbold{1}}_{\{r^\vv(u)=i\}}-
\nu _i) ({\mathbbold{1}}_{\{r^\vv(s)=i\}}-\nu_i)\bigr\}
\,\mathrm{d}s
\\
&\lesssim& \vv^\bb+\E\int_{\bar t_0}^{\tt
t_j}
\bigl\{({\mathbbold{1}}_{\{r^\vv(u)=i\}}-\nu_i) ({
\mathbbold{1}}_{\{
r^\vv
(s)=i\}}-\nu_i)\bigr\}\,\mathrm{d} s,
\end{eqnarray*}
where $\tt t_j:=\bar t_{[(t-t_j)/\vv^\bb]-1}$. Recall from
Yin and Zhang \cite{YZ}, Lemma~7.1, page 169, that
%
\begin{equation}
\label{eq26} \bigl|\P^\vv(u,s)-\nu\bigr|\lesssim \biggl(\vv+\exp \biggl(-
\frac {\kk
(u-s)}{\vv} \biggr) \biggr),\qquad u\ge s\ge0,
\end{equation}
where $\P^\vv(t,s):=(p^\vv_{ij}(u,s))_{1\le i,j\le
n}=(\P(r^\vv(u)=j)|r^\vv(s)=i)_{1\le i,j\le n}$, and $\kk>0$ is
determined by the eigenvalues of $\tt\Q$. Thus, for $\F_{\tt
t_j}:=\si\{r^\vv(s)\dvt 0\le s\le\tt t_j\}$, using the basic property
of conditional expectation, we deduce that
\begin{eqnarray*}
\bigl|\E\bigl\{(\mathbbold{1}_{\{r^\vv(u)=i\}}-\nu_i) (
\mathbbold{1}_{\{r^\vv
(s)=i\}}-\nu_i)\bigr\}\bigr| &\le&\E\bigl(\bigl|
\mathbbold{1}_{\{r^\vv(s)=i\}}-\nu_i\bigr|\cdot\bigl|\bigl(\E (
\mathbbold{1}_{\{
r^\vv(u)=i\}}-\nu_i)|\F_{\tt t_j}\bigr)\bigr|\bigr)
\\
&\lesssim& \E\bigl(\bigl|\bigl(\E(\mathbbold{1}_{\{r^\vv(u)=i\}}-\nu_i)|
\F_{\tt
t_j}\bigr)\bigr|\bigr)
\\
&\lesssim& \biggl(\vv+\exp \biggl(-\frac {\kk(u-\tt t_j)}{\vv} \biggr) \biggr)
\\
&\lesssim &\biggl(\vv+\exp \biggl(-\frac {\kk}{\vv^{1-\bb}} \biggr) \biggr)
\\
&\lesssim& \vv,
\end{eqnarray*}
where in the last third step we used the fact \eqref{eq26}, the last
two step is due to $u>\bar t_{[(t-t_j)/\vv^\bb]}$, while the last
one owes to $\exp(-\frac {\kk}{\vv^{1-\bb}})\lesssim \vv$ for
sufficiently
small $\vv\in(0,1)$. Hence,
%
\begin{equation}
\label{eq19} \eta^\vv(t)\lesssim \vv^{\bb+\rr}.
\end{equation}
Note that from \eqref{eq14} and the uniform boundedness of
$|\mathbbold{1}_{\{r^\vv(s)=i\}}-\nu_i|\leq1$,
\begin{eqnarray*}
&& \biggl(\E \biggl|\int_{t_j}^{t_{j+1}}\{
\mathbbold{1}_{\{r^\vv(s)=i\}
}-\nu _i\}\,\mathrm{d} s
\biggr|^{(p(1+\dd))/\dd} \biggr)^{\dd/(p(1+\dd))}
\\
&&\qquad \le \vv^{\rr-(2\dd\rr)/(p(1+\dd))} \biggl(\E \biggl|\int_{t_j}^{t_{j+1}}
\{ \mathbbold{1}_{\{r^\vv(s)=i\}}-\nu_i\}\,\mathrm{d} s
\biggr|^2 \biggr)^{\dd/(p(1+\dd))}.
\end{eqnarray*}
Then
claim \eqref{eq23} follows from \eqref{eq19}.
Putting \eqref{eq23} and \eqref{eq24} into \eqref{eq25}, we arrive
at
\begin{eqnarray*}
\Upsilon_{3i}(t)&\lesssim _T&\sum
_{j=0}^{\lfloor
t/\vv^\rho\rfloor }\mathrm{e}^{-\lambda_1(t-t_j)}
\vv^{\rr
+(\bb
\dd-2\dd\rr)/(p(1+\dd))}\lesssim _T\bigl(\mathrm{e}^{\lambda_1\vv^\rr}-1
\bigr)^{-1}\vv^{\rr+((\bb-\rr
)\dd)/(p(1+\dd))}
\\
&\lesssim _T&\vv^{((\bb- \rr)\dd)/(p(1+\dd))}
\end{eqnarray*}
due to the fact that $\mathrm{e}^{\lambda_1\vv^\rr
}-1=\mathrm{O}(\lambda_1\vv^\rr)$ for
sufficiently small $\vv\in(0,1)$. So, we get
\begin{eqnarray*}
\bigl(\E\bigl|X^\vv(t)-\bar X(t)\bigr|_H^p
\bigr)^{1/p}&\le& C_T\bigl(\vv^{\rho\theta}+
\vv^{((\bb-\rr)\dd)/(p(1+\dd))}\bigr)
\\
& &{}+ \sum_{i=1}^nK_i\int
_0^t\mathrm{e}^{-\lambda
_1(t-s)}\bigl(
\E\bigl|X^\vv(s)-\bar X(s)\bigr|_H^p
\bigr)^{1/p}\,\mathrm{d}s.
\end{eqnarray*}
It follows from the Gronwall inequality that
\[
\bigl(\E\bigl\|X^\vv(s)-\bar X(s)\bigr\|_H^p
\bigr)^{1/p}\lesssim _T\bigl(\vv^{\rho\theta}+
\vv^{((\bb-\rr)\dd
)/(p(1+\dd))}\bigr).
\]
Then the desired assertion holds by noting that
$\rr<(\aa-\rr)/(\aa-\rr+p\theta\aa)$ and choosing appropriate
$\bb\in(\rho,1)$.
\end{pf*}

%
\begin{rem}
By a close inspection of argument of Theorem~\ref{two}, if
$\sum_{i=1}^nK_i<\lambda_1$, we can also derive a long-term error bound
below
\[
\sup_{t\ge0}\bigl(\E\bigl|X^\vv(t)-\bar
X(t)\bigr|_H^p\bigr)^{1/p}\lesssim
\vv^{\rr
\theta
},\qquad p\in(1,\aa),
\]
for any $\aa\in(1,2)$ and sufficiently small $\vv\in(0,1)$, where
$\theta\in(0,1)$ is the constant such that (A3) holds and
$\rr<(\aa-p)/(\aa-p+p\theta\aa)$.
\end{rem}

\begin{rem}
By means of the martingale problem formulation, the weak
convergence of $(X^\vv(t),\bar r^\vv(t))$
for hybrid finite-dimensional systems were obtained in Yin and Zhang
\cite{YZ}, Theorem~7.20, page~204. In the current framework, it only admits a unique
mild solution rather than strong solution so that the
martingale-problem method seems not to be available. However, in
this subsection, we investigate the strong
convergence (in moment-sense) of $\{X^\vv(t)\}_{t\ge0}$ to the
averaging process $\{\bar X(t)\}_{t\ge0}$ defined by \eqref{eq8}
by the semigroup approach. We also
provide a convergence rate in terms of error bounds. Moreover, even for
$\aa=2$, that is, the Wiener noise case, our result still seems to be
new for infinite-dimensional stochastic dynamical systems.
\end{rem}

\subsection{Two-time-scale Markov switching with multiple weakly
irreducible classes}
In this subsection, we proceed to investigate the averaging
principle associated with \eqref{eq1}, where the Markov chain
$r^\vv(t)$ has a large state space
\[
\S:=\S_1\cup\S_2\cup\cdots\cup\S_l
\]
with $\S_i:=\{s_{i1},\ldots,s_{in_i}\}$ and $n:=n_1+n_2+\cdots+n_l$.
Assume that the generator $\Q^\vv:=(q_{ij}^\vv)_{n\times n}$ of
$r^\vv(t)$ admits the form
\[
\Q^\vv:=\frac {1}{\vv}\tt\Q+\hat\Q,
\]
where $\tt\Q:=(\tt q_{ij})_{n\times
n}=\operatorname{diag}(\tt\Q_1,\ldots,\tt\Q_l)$ such that, for each
$k\in\{1,\ldots,l\}$, $\tt\Q_k$ is irreducible and the generator of
some Markov chain taking values in $\S_k$ with the corresponding
stationary distribution
$\mu_k:=(\mu_{k1},\ldots,\mu_{kn_k})\in\R^{1\times n_k}$, and
$\hat\Q:=(\hat q_{ij})_{n\times n}$. Since the transitions within
each group take place at a fast pace, whereas the interactions from
one group to another are relatively infrequently, following the
basic idea in Yin and Zhang \cite{YZ}, we lump the states in
each $\S_k$ into a single state and then define an aggregated
process $\bar r^\vv(\cdot)$ by
\[
\bar r^\vv(t)=k \qquad\mbox{for } r^\vv(t)\in
\S_k
\]
with the associated state space $\bar\S:=\{1,\ldots,l\}$. Let
\[
\bar\Q:=(\bar q_{ij})_{l\times l}=\tt\mu\hat\Q\mathbf I,
\]
where $\tt\mu:=\operatorname{diag}(\mu_1,\ldots,\mu_l)\in\R^{l\times n}$ and
${\mathbf I}:=\operatorname{diag}({\mathbf I}_{n_1},\ldots,{\mathbf I}_{n_l})$ with ${\mathbf
I}_{n_k}:=(1,\ldots,1)^T\in\R^{n_k\times1},k=1,\ldots,l$. Recall from
Yin and Zhang \cite{YZ}, Theorem~7.4, page~172, that $\bar r^\vv(\cdot)$
converges
weakly to the continuous-time Markov chain $\bar r(\cdot)$ with the
state space $\bar\S$ and the generator $\bar\Q$ as $\vv\to0$,
although generally $\bar r^\vv(t)$ need not be a Markov chain. Our
main result in this subsection is stated as follows.

\begin{thmm}\label{th}
Let \textup{(A1)--(A3)} hold and suppose further that $x\in\mathcal
{D}((-A)^\theta)$. Then
%
\begin{equation}
\lim_{\vv\to0}\E\bigl|X^\vv(t)-\bar X(t)\bigr|_H^p=0,\qquad
t\in[0,T] \mbox{ and } p\in(1,\aa),
\end{equation}
where $\bar X(t)$ satisfies in the mild sense the following
averaging equation
%
\begin{equation}
\label{i1} \mathrm{d}\bar X(t)=\bigl\{A\bar X(t)+\bar b\bigl(\bar X(t),\bar
r(t)\bigr)\bigr\}\,\mathrm{d}t+\mathrm{d} L(t), \qquad\bar X(0)=x, \bar r(0)=\bar
r_0
\end{equation}
with $\bar b(y,i):=\sum_{j=1}^{n_i}\mu_{ij}b(y,s_{ij})$.
\end{thmm}

\begin{pf}
We only give an outline of the proof since it is very similar to
that of Theorem~\ref{two}. By (A1)--(A3), for any $p\in(1,\aa)$, we
deduce that
%
\begin{equation}
\label{i6} \sup_{0\le t\le T}\E\bigl|\bar X(t)\bigr|_H^p<
\infty.
\end{equation}
It is easy to see from (A1) and (A2) that
\begin{eqnarray*}
&&\bigl(\E\bigl|X^\vv(t)-\bar X(t)\bigr|_H^p
\bigr)^{1/p}
\\
& &\quad\le \sum_{i=1}^l\sum
_{j=1}^{n_i}K_{s_{ij}}\int_0^t
\mathrm{e}^{-\lambda_1(t-s)}\bigl(\E\bigl|X^\vv (s)-\bar X(s)\bigr|_H^p
\bigr)^{1/p}\,\mathrm{d}s
\\
&&\quad{} +\sum_{i=1}^l\sum
_{j=1}^{n_i}\biggl(\E \biggl|\int_0^t
\mathrm{e}^{(t-s)A}b\bigl(\bar X(s),s_{ij}\bigr)\{
\mathbbold{1}_{\{r^\vv(s)=s_{ij}\}}-\mu_{ij}\mathbbold {1}_{\{\bar
r^\vv(s)=i\}}\}
\,\mathrm{d}s \biggr|_H^p \biggr)^{1/p}
\\
&&\qquad{} +\sum_{i=1}^l\sum
_{j=1}^{n_i}\mu_{ij} \biggl(\E\biggl|\int
_0^t\mathrm{e} ^{(t-s)A}b\bigl(\bar
X(s),s_{ij}\bigr)\{\mathbbold{1}_{\{\bar
r^\vv(s)=i\}}-\mathbbold{1}_{\{\bar
r(s)=i\}}
\}\,\mathrm{d}s\biggr |_H^p \biggr)^{1/p}
\\
&&\quad =: \Phi_1(t)+\Phi_2(t)+\Phi_3(t).
\end{eqnarray*}
By the definition of $\bar r^\vv(\cdot)$, one has
\[
\bigl\{\bar r^\vv(t)=i\bigr\}=\bigl\{r^\vv(t)\in
\S_i\bigr\}. 
\]
Then, in the same way as the proof of \eqref{eq24}, we deduce
from \eqref{i6} and Yin and Zhang \cite{YZ}, Theorem~7.2, page~170, that
%
\begin{equation}
\label{i2} \Phi_2(t)\to0 \qquad\mbox{as } \vv\to0.
\end{equation}
Next, applying the H\"older inequality, we find that
\begin{eqnarray*}
\Phi_3(t)&\le& \sum_{i=1}^l
\sum_{j=1}^{n_i}\mu_{ij}\int
_0^t\bigl\|\mathrm{e}^{(t-s)A}\bigl({\mathbf
1}-\mathrm{e}^{(s-\lfloor
s\rfloor )A}\bigr)\bigr\|\bigl(\E\bigl|b\bigl(\bar X(s),s_{ij}
\bigr)\bigr|_H^p\bigr)^{1/p}\,\mathrm{d}s
\\
& &{}+\sum_{i=1}^l\sum
_{j=1}^{n_i}\mu_{ij}\int
_0^t\bigl\|\mathrm{e}^{(t-\lfloor
s\rfloor )A}\bigr\|\bigl(\E\bigl|b
\bigl(\bar X(s),s_{ij}\bigr)-b\bigl(\bar X\bigl(\lfloor s\rfloor
\bigr),s_{ij}\bigr)\bigr|_H^p\bigr)^{1/p}
\,\mathrm{d}s
\\
&&{} +\sum_{i=1}^l\sum
_{j=1}^{n_i}\mu_{ij}\sum
_{k=0}^{\lfloor
t/\vv^\rho\rfloor } \bigl(\E\bigl|\mathrm{e}^{(t-t_k)A}b\bigl(
\bar X(t_k),s_{ij}\bigr)\bigr|_H^{p(1+\dd
)}
\bigr)^{1/(p(1+\dd))}
\\
&&\hspace*{40pt}{} \times \biggl(\E \biggl|\int_{t_j}^{t_{j+1}}\{
\mathbbold{1}_{\{
\bar
r^\vv(s)=i\}}-\mathbbold{1}_{\{\bar r(s)=i\}}\}\,\mathrm{d} s
\biggr|^{(p(1+\dd))/\dd} \biggr)^{\dd/(p(1+\dd))}
\\
&=:&\Theta_1(t)+\Theta_2(t)+\Theta_3(t),
\end{eqnarray*}
where $\lfloor s\rfloor :=[s/\vv^\rho]\vv^\rho$ for $\rho\in
(0,1)$, and
$t_j:=j\vv^\rho,j=0,\ldots,[ t/\vv^\rho]$, and $t_{[ t/\vv^p]+1}:=t$.
Moreover, carrying out similar arguments to those of \eqref{i3} and
\eqref{i4} and utilizing \eqref{i6} and Lemma~\ref{Le} yields that
%
\begin{equation}
\label{i5} \Theta_1(t)+\Theta_2(t)\lesssim
\vv^{\rr\theta}.
\end{equation}
From (A1) and \eqref{i6}, for sufficiently small $\vv\in(0,1)$, it
is seen that
\[
\sum_{k=0}^{\lfloor t/\vv^\rho\rfloor } \bigl(\E\bigl|
\mathrm{e}^{(t-t_k)A}b\bigl(\bar X(t_k),s_{ij}
\bigr)\bigr|_H^{p(1+\dd)} \bigr)^{1/(p(1+\dd))}\lesssim
\vv^{-\rr}.
\]
On the other hand, by the weak convergence of $\bar r^\vv(\cdot)$ to
$\bar r(\cdot)$ (Yin and Zhang \cite{YZ}, Theorem~7.4, page~172), the Skorohod
representation theorem (Yin and Zhang \cite{YZ}, Theorem~14.5, page~318), and the
dominated convergence theorem, we have
\[
\biggl(\E \biggl|\int_{t_j}^{t_{j+1}}\{\mathbbold{1}_{\{\bar
r^\vv(s)=i\}}-
\mathbbold{1}_{\{\bar r(s)=i\}}\}\,\mathrm{d} s \biggr|^{(p(1+\dd))/\dd}
\biggr)^{\dd/(p(1+\dd))}\lesssim \vv^\rr g(\vv),
\]
where the positive function $g(\cdot)$ such that $g(\vv)\to0$ as
$\vv\to0$. Then we obtain that
\[
\Theta_3(t)\lesssim g(\vv).
\]
Henceforth the desired assertion follows from the Gronwall
inequality.
\end{pf}

\begin{rem}
Unlike the case discussed in the previous subsection, it seems hard
to give a strong convergence rate bound since the details on $\bar
r(\cdot)$ are not enough, however the averaging equation~\eqref{i1}
is explicitly dependent on the Markov chain $\bar r(\cdot)$, which
is quite different from the case investigated in the last
subsection.
\end{rem}

\section{SPDEs with an additional fast-varying process driven by another
cylindrical stable process}\label{fast-jump} In this section, we
work on another two-time-scale system, in which
there is an additional random process that has a fast-varying
component driven by another cylindrical stable process.

For a small parameter $\vv>0$, we consider the following stochastic
fast--slow system
%
\begin{equation}
\label{f6} \mathrm{d}X^\vv(t)=\bigl\{AX^\vv(t)+ b
\bigl(X^\vv(t),Y^\vv(t)\bigr)\bigr\}\,\mathrm{d}t+
\mathrm{d}L(t),\qquad X^\vv(0)=x\in\mathcal{D}\bigl((-A)^{1/2}
\bigr)
\end{equation}
and
%
\begin{equation}
\label{f7} \mathrm{d}Y^\vv(t)=\frac {1}{\vv}\bigl
\{BY^\vv(t)+f\bigl(X^\vv(t),Y^\vv (t)\bigr)
\bigr\}\,\mathrm{d} t+\frac {1}{\vv^{1/\bb}}\,\mathrm{d}Z(t),\qquad Y^\vv(0)=y
\in H.
\end{equation}

Throughout this section, we shall assume that:
\begin{longlist}[(B4)]
\item[(B1)]$A\dvtx \mathcal{D}(A)\subset H\mapsto H$ is a linear
unbounded operator such that (A1) and $B\dvtx \mathcal{D}(B)\subset
H\mapsto H$
is a self-adjoint compact operator on $H$ such that $-B$ has
discrete spectrum
$0<\mu_1<\mu_2<\cdots<\mu_k<\cdots$ and $\lim_{k\to\infty }\mu
_k=\infty $.

\item[(B2)]
$b$ is uniformly bounded and Lipschitzian, that is, there exist $M,
K_1>0$ such that
\[
\sup_{x,y\in H}\bigl|b(x,y)\bigr|\le M,
\]
and, for $x_1,x_2,y_1,y_2\in H$,
\[
\bigl|b(x_1,y_1)-b(x_2,y_2)\bigr|^2
\le K_1\bigl(|x_1-y_1|^2_H+|x_2-y_2|^2_H
\bigr).
\]

\item[(B3)] For any $x,y\in H$ and $h\in H$, there
exist $K_2,K_3>0$ such that $|\nn^{(1)} f(x,y)\cdot h|\le K_2|h|$
and $|\nn^{(2)} f(x,y)\cdot h|\le K_3|h|$, where $\nn^{(1)} f$ and
$\nn^{(2)} f$ denote the G\^ateaux derivative w.r.t. the first
variable and the second variable, respectively.

\item[(B4)] There exists $\theta\in(0,1)$ such that
$\aa\theta\in(0,1)$,
\[
\kk_1:=\sum_{k=1}^\infty
\frac {\bb_k^\aa}{\lambda_k^{1-\aa
\theta}}<\infty \quad \mbox{and}\quad \kk_2:=\sum
_{k=1}^\infty \frac {q_k^\bb}{\mu_k}<\infty.
\]
\end{longlist}

Under (B1)--(B4), both \eqref{f6} and \eqref{f7} are well-posed in
the mild sense. Consider an SPDE associated with the fast variable,
where the slow variable is fixed and equal to $z\in H$,
%
\begin{equation}
\label{d5} \mathrm{d}Y^z(t;y)=\bigl\{BY^z(t;y)+f
\bigl(z,Y^z(t;y)\bigr)\bigr\}\,\mathrm{d}t+\mathrm{d}Z(t),\qquad
Y^z(0;y)=y\in H.
\end{equation}
Under (B1), (B3) and (B4), \eqref{d5} has a unique mild solution
$\{Y^z(t;y)\}_{t\ge0}$. Moreover, as Lemma~\ref{distr.} below
states, \eqref{d5} admits a unique ergodic invariant measure
$\ppi^z(\cdot)\in\mathcal{P}(H)$, the family of all probability
measures on $H$. Our main result in this section is as follows:

\begin{thmm}\label{fast}
Let \textup{(A1)} and \textup{(B1)--(B4)} hold and assume further that
$K_3<\mu_1$. Then,
%
\begin{equation}
\label{e1} \lim_{\vv\to0}\E\bigl|X^\vv(t)-\bar
X(t)\bigr|_H^p=0,\qquad t\in[0,T], p\in(1,\aa),
\end{equation}
where $\bar X(t)$ is the mild solution of the averaging equation
%
\begin{equation}
\label{e2} \mathrm{d}\bar X(t)=\bigl\{A\bar X(t)+\bar b\bigl(\bar X(t)\bigr)
\bigr\}\,\mathrm{d}t+\mathrm{d}L(t), \qquad \bar X(0)=x\in H
\end{equation}
with
%
\begin{equation}
\label{r4} \bar b(z):=\int_Hb(z,u)\ppi^z(
\mathrm{d}u),\qquad z\in H.
\end{equation}
\end{thmm}

To facilitate the proof of Theorem~\ref{distr.},
we shall present
several technical lemmas in this regards and then finish the
corresponding argument.

\begin{lem}\label{L2}
Under the assumptions of Theorem~\ref{fast},
%
\begin{equation}
\label{e3} \sup_{t\ge0}\E\bigl|Y^\vv(t)\bigr|_H^p<
\infty ,\qquad p\in(1,\aa).
\end{equation}
\end{lem}

\begin{pf}
It is easy to see from Priola and Zabczyk \cite{PZ11}, (4.12), (A1),
and (B1), that
%
\begin{equation}
\label{e4} \sup_{t\ge0}\E\bigl|X^\vv(t)\bigr|_H^p<
\infty.
\end{equation}
Let
\[
\bar Z^\vv(t):=\frac {1}{\vv^{1/\bb}}\int_0^t
\mathrm{e}^{
{(t-s)B}/{\vv}}\,\mathrm{d} Z(s).
\]
By Priola and Zabczyk \cite{PZ11}, (4.12) and (B4), one has
%
\begin{equation}
\label{e6} \E\bigl|\bar Z^\vv(t)\bigr|_H^p\le
\vv^{-p/\bb} \Biggl(\sum_{k=1}^\infty
q_k^\bb\int_0^t
\mathrm{e}^{- {\bb\mu_k(t-s)}/{\vv
}}\,\mathrm{d} s \Biggr)^{p/\bb}\le\bigl(
\bb^{-1}\kk_2\bigr)^{p/\bb}.
\end{equation}
In view of (B1), (B3), \eqref{e4}, and \eqref{e6}, we then derive
that
\begin{eqnarray*}
\bigl(\E\bigl| Y^\vv(t)\bigr|_H^p
\bigr)^{1/p}&\le& |y|_H+\vv^{-1}\int
_0^t \bigl\|\mathrm{e}^{ {(t-s)B}/{\vv
}} \bigr\|\bigl(\E\bigl|f
\bigl(X^\vv (s),Y^\vv(s)\bigr)\bigr|_H^p
\bigr)^{1/p}\,\mathrm{d} s+\bigl(\E\bigl|\bar Z^\vv(t)\bigr|_H^p
\bigr)^{1/p}
\\
&\le& |y|_H+\vv^{-1}\int_0^t
\mathrm{e}^{ {-\mu_1(t-s)}/{\vv}}\bigl\{ c\bigl(1+|z|_H\bigr)+K_3\bigl(
\E \bigl|Y^\vv(s)\bigr|_H^p\bigr)^{1/p}\bigr
\}\,\mathrm{d} s
\\
&\le& c\bigl(1+|y|_H+|z|_H\bigr)+\frac {K_3}{\mu_1}\sup
_{t\ge0}\bigl(\E\bigl| Y^\vv(t)\bigr|_H^p
\bigr)^{1/p}.
\end{eqnarray*}
This therefore leads to \eqref{e3} due to $K_3<\mu_1$.
\end{pf}

\begin{lem}\label{distr.}
Assume that the assumptions of Theorem~\ref{fast} hold. Then
\eqref{d5}
admits a unique ergodic invariant measure $\ppi^z(\cdot)\in\mathcal
{P}(H)$ such that
%
\begin{equation}
\label{r3} \bigl|\E b\bigl(Y^z(t;y)\bigr)-\bar b(z)\bigr|_H
\lesssim \mathrm{e}^{-(\mu
_1-K_3)t}\bigl(1+|y|_H+|z|_H\bigr).
\end{equation}
%
\end{lem}

\begin{pf}
We adopt the remote start method to show existence of an invariant
measure for \eqref{d5}. Let $\hat Z(t):=\sum_{k=1}^\infty q_k\hat
Z_k(t)e_k$, where $\{\hat Z_k(t)\}_{k\ge1}$ is an independent copy
of $\{ Z_k(t)\}_{k\ge1}$, and $\{\tt Z(t)\}_{t\ge0}$ be a
double-sided cylindrical $\bb$-stable process defined by
\[
\tt Z(t):= %
\cases{ Z(t), &\quad $t\ge0,$\vspace*{2pt}
\cr
\hat Z(-t),&\quad $t<0$}
\]
with the filtration
\[
\bar\F_t:=\bigcap_{s>t}\bar
\F^0_s,
\]
where $\bar\F^0_s:=\si(\{\tt Z(r_2)-\tt Z(r_1)\dvt -\infty <r_1\le
r_2\le
s,\GG\},\mathscr{N})$ and $\mathscr{N}:=\{A\in\F|\P(A)=0\}$. Next,
consider \eqref{d5}, for arbitrary $s\in(-\infty ,t]$ with $t\in\R$,
%
\begin{equation}
\label{05} \mathrm{d}Y^z(t;s,y)=\bigl\{BY^z(t;s,y)+f
\bigl(z,Y^z(t;s,y)\bigr)\bigr\}\,\mathrm{d}t+\mathrm{d}\tt Z(t),\qquad
Y^z(s;s,y)=y\in H.
\end{equation}
Set $\Gamma^z(t;y):=Y^z(t;-\lambda,y)-Y^z(t;-\gamma,y)$ for
$-\lambda\in(-\gamma,
t]$. By (B1) and (B3), it follows that
\begin{eqnarray*}
\bigl(\E\bigl|\Gamma^z(t;y)\bigr|_H^p
\bigr)^{1/p}\le\mathrm{e}^{-\mu_1(t+\lambda
)}\bigl(\E\bigl|\Gamma^z(-
\lambda ;y\bigr)\bigr|_H^p\bigr)^{1/p}+K_3
\int_{-\lambda}^t\mathrm{e}^{-\mu
_1(t-s)}\bigl(\E\bigl|
\Gamma ^z(s;y)\bigr|_H^p\bigr)^{1/p}\,
\mathrm{d} s.
\end{eqnarray*}
Multiplying $\mathrm{e}^{\mu_1t}$ on both sides leads to
\begin{eqnarray*}
\mathrm{e}^{\mu_1t}\bigl(\E\bigl|\Gamma^z(t;y)\bigr|_H^p
\bigr)^{1/p}\le\mathrm{e}^{-\mu_1\lambda}\bigl(\E\bigl|\Gamma
^z(-\lambda;y)\bigr|_H^p\bigr)^{1/p}+K_3
\int_{-\lambda}^t\mathrm{e}^{\mu_1s}\bigl(\E\bigl|
\Gamma ^z(s;y)\bigr|_H^p\bigr)^{1/p}\,
\mathrm{d} s.
\end{eqnarray*}
Thus we get from the Gronwall inequality that
%
\begin{equation}
\label{01} %
\bigl(\E\bigl|\Gamma^z(t;y)\bigr|_H^p
\bigr)^{1/p}\le\mathrm{e}^{-(\mu
_1-K_3)(t+\lambda)}\bigl(\E\bigl|\Gamma
^z(-\lambda;y)\bigr|_H^p\bigr)^{1/p}.
\end{equation}
Moreover, carrying out an argument of Lemma~\ref{L2}, we have
%
\begin{equation}
\label{02} \sup_{t\ge s}\bigl(\E\bigl|Y^z(t;s,y)\bigr|_H^p
\bigr)^{1/p}\lesssim 1+|y|_H+|z|_H,\qquad s\in \R.
\end{equation}
For $t=0$ and $-\lambda\in(-\mu,0]$, we deduce from \eqref{01} and
\eqref{02} that
\[
\bigl(\E\bigl|Y^z(0;-\lambda,y)-Y^z(0;-
\gamma,y)\bigr|_H^p\bigr)^{1/p}\lesssim
\bigl(1+|y|_H+|z|_H\bigr)\mathrm{e}^{-(\mu
_1-K_3)\lambda}.
\]
From the estimate above, we conclude that $\{Y^z(0;-t,y)\}_{t\ge0}$
is a Cauchy sequence in $L^p(\OO;H)$, and therefore it is convergent
to a random variable $\eta_z(y)\in L^p(\OO;H)$, which is independent
of $y\in H$, and denoted by $\eta_z\in L^p(\OO;H)$. Then, following
a standard procedure (see, e.g., Pr\'{e}v\^ot and R\"ockner \cite{PR07}, pages~109--110), we deduce
that $\mathscr{L}(\eta_z)=:\ppi^z(\cdot)$ is an invariant measure of~\eqref{d5}.

Next, following an argument of \eqref{01}, we obtain that
%
\begin{equation}
\label{r1} \bigl(\E\bigl|Y^z(t;y_1)-Y^z(t;y_2)\bigr|_H^p
\bigr)^{1/p}\le\mathrm{e}^{-(\mu_1-
K_3)t}|y_1-y_2|_H.
\end{equation}
This, together with \eqref{02}, implies that
%
\begin{equation}
\label{06} \E\bigl|Y^z(t;y)\bigr|_H^p\le
\mathrm{e}^{-p(\mu_1- K_3)t}|y|_H^p+c\bigl(1+|z|_H^p
\bigr).
\end{equation}
Furthermore, by virtue of \eqref{06} and using a stationary solution
$Y^z(t,y)$ with invariant law $\ppi^z(\cdot)$, we obtain that
%
\begin{equation}
\label{03} \E^z|y|_H^p =
\E^z\bigl|Y^z(t,y)\bigr|_H^p\le
\mathrm{e}^{-p(\mu_1-
K_3)t}\E^z|y|_H^p+c
\bigl(1+|z|_H^p\bigr),\qquad t\ge0,
\end{equation}
where $\E^z$ is the mathematical expectation operator w.r.t.
$\ppi^z(\cdot)$. \eqref{03} further gives that
%
\begin{equation}
\label{r2} \ppi^z\bigl(|\cdot|_H^p
\bigr)\lesssim 1+|z|_H^p.
\end{equation}
Consequently, \eqref{r1} and \eqref{r2} yield the uniqueness of
invariant measure. Indeed,
if $\tt\ppi^z(\cdot)\in\mathcal{P}(H)$ is also an invariant
measure, for
any $\psi\in C_b^2(H;\R)$, by the
invariance of $ \ppi^z(\cdot)$ and $\tt\ppi^z(\cdot)$, we deduce from
\eqref{r1} and \eqref{r2} that
\begin{eqnarray*}
\bigl|\ppi^z(\psi)-\tt\ppi^z(\psi)\bigr| 
&
\le& c\mathrm{e}^{-(\mu_1- K_3)t}\bigl\{\ppi^z\bigl(|
\cdot|_H\bigr)+\tt \ppi\bigl(|\cdot|_H\bigr)\bigr\}
\\
&\le& c\mathrm{e}^{-(\mu_1-K_3)t}\bigl\{1+|z|_H\bigr\}
\to0 \qquad\mbox{as } t\uparrow\infty.
\end{eqnarray*}
That is, for any $\psi\in C_b^2(H;\R)$, $\ppi^z(\psi)=\tt
\ppi^z(\psi)$, which shows that $\ppi\equiv\tt\ppi$ due to
Ikeda and Watanabe \cite{IW}, Proposition~2.2, page~3.

Finally, \eqref{r3} follows by
noting from the invariance of $\ppi^z(\cdot)$, \eqref{r2} and the
Lipschitz property of $b$.
\end{pf}

Applying the Lipschitzian property of $b$, the ergodic property of
invariant measure $\ppi^z(\cdot)\in\mathcal{P}(H)$ due to Lemma~\ref{distr.} and the uniform boundedness of the directional
derivative $\nn_hY^z(t;y)$ with respect to $z\in H$ along the
direction $h\in H$, and adopting a similar argument in
Cerrai and Freidlin \cite{CF}, (5.4), we deduce that $\bar b$ is
Lipschitzian, which is
stated as the following corollary for citation convenience.

\begin{cor}\label{Lip}
Under the assumptions of Theorem~\ref{fast}, $\bar b\dvtx H\to H$ is
Lipschitzian.
\end{cor}

To reveal the error analysis between the slow component
$\{X^\vv(t)\}_{t\ge0}$ and the averaging process $\{\bar
X(t)\}_{t\ge0}$, determined by \eqref{e2}, we further need to define
the following two auxiliary processes:
%
\begin{equation}
\label{g3} \tt Y^\vv(t):=\mathrm{e}^{ {tB}/{\vv}}y+
\frac {1}{\vv}\int_0^t\mathrm{e}^{ {(t-s)B}/{\vv
}}f
\bigl(X^\vv\bigl(\lfloor s/\dd\rfloor \dd\bigr),\tt
Y^\vv(s)\bigr)\,\mathrm{d} s+\frac {1}{\vv^{1/\bb}}\int
_0^t\mathrm{e}^{
{(t-s)B}/{\vv}}\,\mathrm{d}Z(s)
\end{equation}
and
%
\begin{equation}
\label{g2} \tt X^\vv(t):=\mathrm{e}^{tA}x+\int
_0^t\mathrm{e}^{(t-s)A}b
\bigl(X^\vv\bigl(\lfloor s/\dd\rfloor \dd\bigr),\tt
Y^\vv(s)\bigr)\,\mathrm{d}s+\int_0^t
\mathrm{e}^{(t-s)A}\,\mathrm{d}L(s),
\end{equation}
where $\dd\in(\vv,1)$ is some constant to be chosen.

\begin{lem}\label{L1}
Assume that the assumptions of Theorem~\ref{fast} hold. Then,
for any $p\in(1,\aa)$,
%
\begin{equation}
\label{h2} \int_0^T\bigl(
\E\bigl|Y^\vv(s)-\tt Y^\vv(s) \bigr|_H^p
\bigr)^{1/p}\,\mathrm{d}s \lesssim \frac {\vv}{\dd}+\vv
\dd^{-(1-\theta)}\mathrm{e}^{ {K_3\dd
}/{\vv}}
\end{equation}
and
%
\begin{equation}
\label{h3} \int_0^T\bigl(
\E\bigl|X^\vv(s)-\tt X^\vv(s)\bigr |_H^p
\bigr)^{1/p}\,\mathrm{d}s \lesssim \dd^\theta+
\frac {\vv}{\dd}+\vv\dd^{-(1-\theta)}\mathrm{e}^{ {
K_3\dd}/{\vv}},
\end{equation}
where $\theta\in(0,1)$ is the constant such that \textup{(B4)}.
\end{lem}

\begin{pf}
For notation simplicity, we set $\Lambda^\vv(t):=Y^\vv(t)-\tt
Y^\vv(t)$. By Lemma~\ref{Le},
for any $t\in[0,T]$
it follows from (B1) and (B2) that
\begin{eqnarray*}
\bigl(\E\bigl|X^\vv(t)-\tt X^\vv(t) \bigr|_H^p
\bigr)^{1/p}&\le&\int_0^t\bigl(\E\bigl|b
\bigl(X^\vv(s),Y^\vv(s)\bigr)-b\bigl(X^\vv
\bigl(\lfloor s/\dd\rfloor \dd\bigr),\tt Y^\vv(s)\bigr)\bigr|_H^p
\bigr)^{1/p}\,\mathrm{d}s
\\
&\lesssim&\int_0^t\bigl(\E\bigl|X^\vv(s)-X^\vv
\bigl(\lfloor s/\dd\rfloor \dd\bigr)\bigr|_H^p
\bigr)^{1/p}\,\mathrm{d}s+\int_0^t
\bigl(\E\bigl|\Lambda^\vv (s)\bigr|_H^p
\bigr)^{1/p}\,\mathrm{d} s
\\
&\lesssim _T&\dd^\theta+\int_0^t
\bigl(\E\bigl|\Lambda^\vv (s)\bigr|_H^p
\bigr)^{1/p}\,\mathrm{d}s.
\end{eqnarray*}
Therefore, to complete the proof of Lemma~\ref{L1}, it is sufficient
to show \eqref{h2}. Carrying out similar arguments to those of
\eqref{e3} and \eqref{e4}, we also deduce that
%
\begin{equation}
\label{h1} \sup_{t\ge0}\E\bigl|\tt X^\vv(t)\bigr|_H^p
\vee\sup_{t\ge0}\E\bigl|\tt Y^\vv(t)\bigr|_H^p<
\infty.
\end{equation}
For any $t\in[0,T]$, there exists an integer $k\ge0$ such that
$t\in[k\dd,(k+1)\dd)$. From (B3) and Lemma~\ref{Le}, we derive that
\begin{eqnarray*}
&&\bigl(\E\bigl|\Lambda^\vv(t) \bigr|_H^p
\bigr)^{1/p}
\\
&&\quad\le\mathrm{e}^{ {-\mu_1(t-k\dd)}/{\vv}}\bigl(\E \bigl|\Lambda^\vv(k\dd
)\bigr|_H^p\bigr)^{1/p}\\
&&\qquad{}+\frac {1}{\vv}\int
_{k\dd}^t\mathrm{e}^{
{-\mu
_1(t-s)}/{\vv}}\bigl(\E \bigl|f
\bigl(X^\vv(s),Y^\vv(s)\bigr)-f\bigl(X^\vv(k
\dd),\tt Y^\vv (s)\bigr)\bigr|_H^p
\bigr)^{1/p}\,\mathrm{d}s
\\
&&\quad\le \mathrm{e}^{ {-\mu_1(t-k\dd)}/{\vv}}\bigl(\E\bigl|\Lambda^\vv (k\dd
)\bigr|_H^p\bigr)^{1/p}\\
&&\qquad{}+\frac
{1}{\vv}\int
_{k\dd}^t\mathrm{e}^{ {-\mu_1(t-s)}/{\vv}}\bigl\{
K_2\bigl(\E\bigl|X^\vv(s)-X^\vv (k
\dd)\bigr|_H^p\bigr)^{1/p}+K_3\bigl(
\E\bigl|\Lambda^\vv(s)\bigr|_H^p\bigr)^{1/p}
\bigr\}\,\mathrm{d} s.
\end{eqnarray*}
This, together with the combined use of \eqref{e3} and \eqref{h1},
yields that
\begin{eqnarray*}
&&\mathrm{e}^{ {\mu_1t}/{\vv}}\bigl(\E\bigl|\Lambda^\vv(t)
\bigr|_H^p\bigr)^{1/p}
\\
&&\quad\le c\mathrm{e}^{ {\mu_1k\dd}/{\vv}}+\frac {c}{\vv}\int_{k\dd
}^t
\mathrm{e}^{ {\mu_1s}/{\vv
}}\bigl(\E\bigl|X^\vv(s)-X^\vv(k
\dd)\bigr|_H^p\bigr)^{1/p}+\frac {
K_3}{\vv}\int
_{k\dd}^t\mathrm{e}^{ {\mu_1s}/{\vv}}\bigl(\E \bigl|
\Lambda^\vv (s)\bigr|_H^p\bigr)^{1/p}\,
\mathrm{d} s.
\end{eqnarray*}
Then, applying the Gronwall inequality and using Lemma~\ref{Le}, we
obtain that
\begin{eqnarray*}
&&\bigl(\E\bigl|\Lambda^\vv(t) \bigr|_H^p
\bigr)^{1/p}
\\
&&\quad\le c\mathrm{e}^{- {(\mu_1- K_3)(t-k\dd)}/{\vv}}\\
&&\qquad{}+\frac {
K_3}{\vv} \int_{k\dd}^t
\mathrm{e}^{ {(-(\mu_1- K_3)t-
K_3k\dd+\mu_1s)}/{\vv}}\bigl(\E\bigl|X^\vv(s)-X^\vv(k
\dd)\bigr|_H^p\bigr)^{1/p}\,\mathrm{d} s
\\
&&\quad\le c\mathrm{e}^{- {(\mu_1-K_3)(t-k\dd)}/{\vv}}-\frac {c
K_3\dd^\theta}{\lambda_1} \mathrm{e}^{ {-(\mu
_1-K_3)(t-k\dd)}/{\vv}} +
\frac {c K_3\dd^\theta}{\mu_1} \mathrm{e}^{ {K_3(t-k\dd)}/{\vv}}
\\
&&\quad\lesssim \mathrm{e}^{- {(\mu_1-K_3)(t-k\dd)}/{\vv
}}+\frac { K_3\dd
^\theta}{\mu_1} \mathrm{e}^{ {K_3(t-k\dd)}/{\vv}}.
\end{eqnarray*}
Integrating from $k\dd$ to $(k+1)\dd$ with respect to the variable
$t$ in the above leads to
\begin{eqnarray*}
\int_{k\dd}^{(k+1)\dd}\bigl(\E\bigl|\Lambda^\vv(t)\bigr|_H^p
\bigr)^{1/p}\,\mathrm{d} t&\lesssim &\int_{k\dd}^{(k+1)\dd}
\biggl\{\mathrm{e}^{-
{(\mu_1-
K_3)(t-k\dd)}/{\vv}}+\frac { K_3\dd^{1/2}}{\lambda_1} \mathrm{e}^{ {
K_3(t-k\dd)}/{\vv}}
\biggr\}\,\mathrm{d} t\\
&\lesssim& \vv+\vv\dd^\theta\mathrm{e}
^{ {K_3\dd}/{\vv}}.
\end{eqnarray*}
Thus, \eqref{h2} follows.
\end{pf}

\begin{rem}
Br\'{e}hier \cite{B12}, Lemma~3.1, confined Lemma~\ref{L1} on
the case $p=1$, which is not sufficient for our purposes, and the
techniques used there does not work for our model. On the other
hand, for finite-dimensional jump--diffusion processes,
Givon \cite{G07}, Lemma~2.4, gives a similar estimate making use of the
It\^o formula, which is unavailable for our framework since the
noise process does not admits second moments.
\end{rem}

With the previous lemmas at hand, we now can complete the
proof of Theorem~\ref{fast}.

\begin{pf*}{Proof of Theorem~\ref{fast}} The proof is inspired by
Khasminskii \cite{K68}.
According to (B2), Lemmas~\ref{Le} and~\ref{L1}, it then
follows that
\begin{eqnarray*}
\bigl(\E\bigl|X^\vv(t)-\bar X(t) \bigr|_H^p
\bigr)^{1/p} &\lesssim& \bigl(\E\bigl|\tt X^\vv(t)-\bar X(t)
\bigr|_H^p\bigr)^{1/p}+\int_0^t
\bigl(\E\bigl|X^\vv(s)-X^\vv\bigl(\lfloor s/\dd\rfloor \dd
\bigr)\bigr|_H^p\bigr)^{1/p}\,\mathrm{d}s
\\
&&{} +\int_0^t\bigl(\E\bigl|Y^\vv(s)-\tt
Y^\vv(s)\bigr|_H^p\bigr)^{1/p}\,
\mathrm{d}s
\\
&\lesssim& \dd^\theta+\frac {\vv}{\dd}+\vv\dd^{-(1-\theta
)}
\mathrm{e}^{ {
K_3\dd}/{\vv}}+\bigl(\E\bigl|\tt X^\vv(t)-\bar X(t)
\bigr|_H^p\bigr)^{1/p}.
\end{eqnarray*}
Therefore, to get the desired assertion, it is sufficient to show
that\label{page20}
%
\begin{equation}
\label{w1} \bigl(\E\bigl|\tt X^\vv(t)-\bar X(t) \bigr|_H^p
\bigr)^{1/p}\lesssim \dd^\theta+\frac {\vv}{\dd}+\sqrt{
\frac
{\vv}{\dd
}}+\vv\dd ^{-(1-\theta)}\mathrm{e}^{ {
K_3\dd}/{\vv}}.
\end{equation}
By the Lipschitz property of $\bar b$ due to Corollary~\ref{Lip},
Lemmas \ref{Le} and \ref{L1}, we
deduce that
%
\begin{eqnarray}
\label{t5} %
&&\bigl(\E\bigl|\tt X^\vv(t)-\bar X(t)
\bigr|_H^p\bigr)^{1/p}\nonumber
\\
&&\quad\le \biggl(\E \biggl|\int_0^t\mathrm{e}^{(t-s)A}
\bigl\{ b\bigl(X^\vv\bigl(\lfloor s/\dd\rfloor \dd\bigr),\tt Y(s)
\bigr)-\bar b\bigl(X^\vv\bigl(\lfloor s/\dd\rfloor \dd\bigr)\bigr)
\bigr\}\,\mathrm{d}s \biggr|_H^p \biggr)^{1/p}\nonumber
\\
&&\qquad{} +\int_0^t\bigl(\E\bigl|\bar b
\bigl(X^\vv\bigl(\lfloor s/\dd\rfloor \dd\bigr)\bigr)-\bar b\bigl(
X^\vv(s)\bigr)\bigr|_H^p\bigr)^{1/p}\,
\mathrm{d}s
\nonumber
\\[-8pt]
\\[-8pt]
\nonumber
&&\qquad{} +\int_0^t\bigl(\E\bigl|\bar b\bigl(
X^\vv(s)\bigr)-\bar b\bigl(\tt X^\vv(s)
\bigr)\bigr|_H^p\bigr)^{1/p}\,\mathrm{d}s+\int
_0^t\bigl(\E\bigl|\bar b\bigl(\tt X^\vv
(s)\bigr)-\bar b\bigl(\bar X(s)\bigr)\bigr|_H^p
\bigr)^{1/p}\,\mathrm{d}s
\\
&&\quad\lesssim \dd^\theta+\frac {\vv}{\dd}+\vv\dd^{-(1-\theta
)}
\mathrm{e}^{ {
K_3\dd}/{\vv}}+\int_0^t\bigl(\E\bigl|\tt X^\vv(s)-\bar X(s)\bigr|_H^p\bigr)^{1/p}
\,\mathrm{d}s\nonumber
\\
&&\qquad{} + \biggl(\E \biggl|\int_0^t\mathrm{e}^{(t-s)A}
\bigl\{b\bigl(X^\vv\bigl(\lfloor s/\dd\rfloor \dd\bigr),\tt Y(s)
\bigr)-\bar b\bigl(X^\vv\bigl(\lfloor s/\dd\rfloor \dd\bigr)\bigr)
\bigr\}\,\mathrm{d}s \biggr|_H^p \biggr)^{1/p}.\nonumber
\end{eqnarray}
Furthermore, noting that
\[
\biggl|\int_0^th(s)\,\mathrm{d}s
\biggr|^2_H=2\int_0^t\int
_s^t \bigl\langle h(r),h(s)\bigr\rangle
_H\,\mathrm{d} r\,\mathrm{d}s
\]
for a locally integrable function $h\dvtx [0,\infty )\mapsto H$, we obtain
from Jensen's inequality that
%
\begin{eqnarray}
\label{t4} %
& &\biggl(\E \biggl|\int_0^t
\mathrm{e}^{(t-s)A}\bigl\{b\bigl(X^\vv\bigl(\lfloor s/\dd
\rfloor \dd\bigr),\tt Y(s)\bigr)-\bar b\bigl(X^\vv\bigl(\lfloor s/\dd
\rfloor \dd\bigr)\bigr)\bigr\}\,\mathrm{d}s \biggr|_H^p
\biggr)^{1/p}\nonumber
\\
&&\quad\le\sum_{k=0}^{\lfloor
t/\dd\rfloor }
\biggl(\E \biggl|\int_{k\dd}^{(k+1)\dd}\mathrm{e}^{(t-s)A}
\bigl\{b\bigl(X^\vv(k\dd ),\tt Y(s)\bigr)-\bar b\bigl(X^\vv(k
\dd)\bigr)\bigr\}\,\mathrm{d}s\biggr |_H^p
\biggr)^{1/p}
\\
&&\quad\lesssim \vv\sum_{k=0}^{\lfloor
t/\dd\rfloor } \biggl(\int
_0^{\dd/\vv}\int_s^{\dd/\vv}
\mathcal {J}_k(r,s)\,\mathrm{d} r\,\mathrm{d}s \biggr)^{1/2},\nonumber
\end{eqnarray}
where $t:=(\lfloor t/\dd\rfloor +1)\dd$ and
\begin{eqnarray*}
\mathcal{J}_k(r,s)&:=&\E\bigl\langle \mathrm{e}^{(t-(k\dd+r\vv
))A}
\bigl(b\bigl(X^\vv (k\dd),\tt Y(r\vv+k\dd)\bigr)-\bar b
\bigl(X^\vv(k\dd)\bigr)\bigr),
\\
& &\hspace*{10pt}{}\mathrm{e}^{(t-(k\dd+s\vv))A}\bigl(b\bigl(X^\vv(k\dd),\tt Y(s\vv +k\dd)
\bigr)-\bar b\bigl(X^\vv(k\dd)\bigr)\bigr)\bigr\rangle _H.
\end{eqnarray*}
For any $s\in(0,\dd)$, observe from \eqref{g3} that
%
\begin{eqnarray}
\label{w7} %
\tt Y^\vv(s+k\dd) &=&\mathrm{e}^{ {sB}/{\vv}}
\tt Y^\vv(k\dd)+\frac {1}{\vv}\int_0^{s}
\mathrm{e}^{
{(s-u)B}/{\vv
}}f\bigl(X^\vv(k\dd),\tt Y^\vv(k
\dd+u)\bigr)\,\mathrm{d}u
\nonumber
\\[-8pt]
\\[-8pt]
\nonumber
&&{} +\frac {1}{\vv^{1/\bb}}\int_0^{s}
\mathrm{e}^{
{(s-u)B}/{\vv}}\,\mathrm{d}Z_1(u), %
\end{eqnarray}
where $Z_1(\cdot):=Z(\cdot+k\dd)-Z(k\dd)$ with filtration
$\F_{\cdot+k\dd}$, which is again a cylindrical $\bb$-stable
process. Let
\[
Z_2(t):=\sum_{k=1}^\infty
q_k \bar Z_{k}(t)e_k,
\]
where $\{\bar Z_{k}(t)\}_{k\ge1}$ is a sequence of i.i.d.
$\R$-valued symmetric $\bb$-stable L\'{e}vy processes defined on the
filtered probability space $(\OO, \F, {\{\F_t}\}_{t\ge0}, \P)$ such
that $ \{Z_2(t)\}_{t\ge0}$ is independent of $\{L(t)\}_{\ge0}$ and
$\{Z(t)\}_{t\ge0}$, respectively. For fixed $X^\vv(k\dd)$ and the
starting point $\tt Y^\vv(k\dd)$, define the process
$Y_s^{X^\vv(k\dd),\tt Y^\vv(k\dd)}$ by
%
\begin{eqnarray}
\label{w8} %
Y_{s/\vv}^{X^\vv(k\dd),\tt Y^\vv(k\dd)}&:=&
\mathrm{e}^{
{sB}/{\vv} }\tt Y^\vv(k\dd)+\int_0^{s/\vv}
\mathrm{e}^{( {s}/{\vv
}-u)B}f \bigl(X^\vv(k\dd ),Y_u^{X^\vv(k\dd),\tt
Y^\vv(k\dd)}
\bigr)\,\mathrm{d}u
\nonumber
\\[-8pt]
\\[-8pt]
\nonumber
&&{} +\int_0^{s/\vv}\mathrm{e}^{( {s}/{\vv}-u)B}\,
\mathrm{d}Z_2(u). %
\end{eqnarray}
A simple calculation gives that
%
\begin{eqnarray}
\label{w8-a} %
Y_{s/\vv}^{X^\vv(k\dd),\tt Y^\vv(k\dd)}&=&
\mathrm{e}^{
{sB}/{\vv} }\tt Y^\vv(k\dd)+\frac {1}{\vv}\int
_0^{s}\mathrm{e}^{
{(s-u)B}/{\vv
}}f
\bigl(X^\vv(k\dd ),Y_{u/\vv}^{X^\vv(k\dd),\tt
Y^\vv(k\dd)}\bigr)\,
\mathrm{d}u
\nonumber
\\[-8pt]
\\[-8pt]
\nonumber
&&{} +\frac {1}{\vv^{1/\bb}}\int_0^{s}
\mathrm{e}^{
{(s-u)B}/{\vv}}\,\mathrm{d} Z_3(u),\qquad s\in(0,\dd),
\end{eqnarray}
where $ Z_3(\cdot):=\vv^{1/\bb}Z_2(\cdot/\vv)$. By the self-similar
property of stable L\'{e}vy processes (Applebaum~\cite{A09}, page~51), we conclude
from \eqref{w7} and \eqref{w8} that
%
\begin{equation}
\label{t1} \mathscr{L}\bigl(\tt Y^\vv(s+k\dd)\bigr)=\mathscr{L}
\bigl(Y_{s/\vv}^{X^\vv(k\dd),\tt
Y^\vv(k\dd)} \bigr), \qquad s\in(0,\dd).
\end{equation}
This
further implies from \eqref{h1} that
%
\begin{equation}
\label{t2} \sup_{s\in[0,\dd]}\E\bigl |Y_s^{X^\vv(k\dd),\tt
Y^\vv(k\dd)}
\bigr|_H^p<\infty.
\end{equation}
Let
\[
\F_s:=\si\bigl\{Y_u^{X^\vv(k\dd),\tt Y^\vv(k\dd)},u\le s\bigr\}.
\]
Then $X^\vv(k\dd)\in\F_s$. By the property of conditional
expectation (Applebaum \cite{A09}, Lemma~1.1.9), and the boundedness of
$b$ due
to (B2), for $r>s$ we obtain from \eqref{t1} that
\begin{eqnarray*}
\mathcal{J}_k(r,s) &=&\E \bigl\langle \mathrm{e}^{(t-(k\dd+s\vv))A}
\bigl(b \bigl(X^\vv (k\dd),Y_s^{X^\vv(k\dd
),\tt
Y^\vv(k\dd)} \bigr)-
\bar b\bigl(X^\vv(k\dd)\bigr) \bigr)
\\
&&\hspace*{7pt}{} \times\mathrm{e}^{(t-(k\dd+r\vv))A}\bigl( \E \bigl(b \bigl(X^\vv(k\dd
),Y_r^{X^\vv(k\dd),\tt
Y^\vv(k\dd)} \bigr)-\bar b\bigl(X^\vv(k\dd)
\bigr) \bigr) |\F_s \bigr) \bigr\rangle _H
\\
&=&\E \bigl\langle \mathrm{e}^{(t-(k\dd+s\vv))A} \bigl(b \bigl(X^\vv (k
\dd),Y_s^{X^\vv(k\dd
),\tt
Y^\vv(k\dd)}\bigr)-\bar b\bigl(X^\vv(k\dd)
\bigr) \bigr)
\\
&&\hspace*{7pt}{} \times\mathrm{e}^{(t-(k\dd+r\vv))A} \bigl(\E  \bigl(b \bigl(z_1,Y_{r-s}^{X^\vv(k\dd),\tt
Y^\vv(k\dd)}+z_2
\bigr)-\bar b(z_1) \bigr) \bigr) |_{z_2=Y_s^{X^\vv(k\dd),\tt Y^\vv(k\dd
)}}^{z_1=X^\vv
(k\dd)}
\bigr\rangle _H
\\
&\le& \bigl(\E \bigl| b \bigl(z_1,Y_s^{X^\vv(k\dd),\tt
Y^\vv(k\dd)}
\bigr)-\bar b\bigl(z_1(\xi)\bigr)  \bigr|^2_H
\bigr)^{1/2}
\\
&&{} \times \bigl(\E \bigl|\bigl(\E  \bigl(b \bigl(z_1,Y_{r-s}^{X^\vv
(k\dd),\tt
Y^\vv(k\dd)}+z_2
\bigr)-\bar b\bigl(z_1(\xi) \bigr) \bigr) \bigr)
|_{z_2=Y_s^{X^\vv(k\dd),\tt Y^\vv
(k\dd
)}}^{z_1(\xi)=X^\vv(k\dd)} \bigr|^2_H
\bigr)^{1/2}
\\
&\lesssim& \E \bigl(\bigl | \bigl(\E  \bigl(b \bigl(z_1,Y_{r-s}^{X^\vv
(k\dd),\tt
Y^\vv(k\dd)}+z_2
\bigr)-\bar b(z_1) \bigr) \bigr) |_{z_2=Y_s^{X^\vv(k\dd),\tt
Y^\vv(k\dd)}}^{z_1=X^\vv(k\dd)}
\bigr|_H \bigr),
\end{eqnarray*}
where in the last step we have used the boundedness of $b$ due to
(B2). The previous estimation, combining Lemma~\ref{distr.} with
\eqref{e4} and \eqref{t2}, yields that
%
\begin{eqnarray}
\label{t3} %
\mathcal{J}_k(r,s) &\lesssim&
\mathrm{e}^{-(\mu_1-
K_3)(r-s)}\E \bigl(1+\bigl|X^\vv(k\dd)\bigr|_H+
\bigl|Y_s^{X^\vv(k\dd),\tt
Y^\vv(k\dd)}\bigr |_H \bigr)
\nonumber
\\[-8pt]
\\[-8pt]
\nonumber
&\lesssim& \mathrm{e}^{-(\mu_1- K_3)(r-s)}. %
\end{eqnarray}
Thus \eqref{w1} follows by putting \eqref{t3} into \eqref{t4} and
applying the Gronwall inequality in \eqref{t5}. Hence, we obtain that
\[
\bigl(\E\bigl|X^\vv(t)-\bar X(t) \bigr|_H^p
\bigr)^{1/p} \lesssim \dd^\theta+\frac {\vv}{\dd}+\sqrt{
\frac {\vv}{\dd
}}+\vv\dd ^{-(1-\theta)}\mathrm{e} ^{ {
K_3\dd}/{\vv}}.
\]
Letting $\dd:=\vv(-\ln\vv)^{ {1}/{2}}$ and then taking $\vv
\to0$
yields the desired assertion, as required. 
\end{pf*}

\begin{rem}
In this section, we show an averaging result for a class of
two-time-scale SPDEs driven by cylindrical stable noises in the
abstract setting. Therefore, stochastic evolution equations of
parabolic type with slow and fast time scales fit into our
framework.
\end{rem}

\begin{rem}
If $\aa=2$ and $\bb=2$ in Theorem~\ref{fast}, which corresponds
to the cylindrical Wiener noises, by reexamining the argument of
Theorem~\ref{fast}, the boundedness of $ b$ can be removed by
imposing, for example,
\[
\bigl| f(x,y)\bigr|_H\le c_1+c_2|y|,\qquad x,y\in H
\]
for some
appropriate constants $c_1,c_2>0$, that is, $f$ is uniformly bounded
w.r.t. the first variable. Moreover, by a close inspection of
argument of Theorem~\ref{fast}, the boundedness of second moment of
$X^\vv$
plays an important role in error analysis.
However, for
the case $\aa,\bb\in(1,2)$, $X^\vv(\cdot)$ only has the $p$th
moment with $p\in(1,\aa)$. Therefore, for the technical reason, it
seems hard to show Theorem~\ref{fast} without the uniform
boundedness of the nonlinearity. However, for the weak convergence
(e.g., convergence in probability) of averaging principle for systems
\eqref{f6} and \eqref{f7}, the boundedness of the nonlinearity can
be removed. Such result will be reported in our forthcoming
paper.
\end{rem}

\begin{rem}
In this section, we aim to obtaining averaging principles for a
class of SPDEs driven by $\aa$-stable noise with $\aa\in(1,2]$.
However, for the case $\aa\in(0,1)$ the method of this paper does
not work. For such a case, it is necessary
to
find new
approaches for the investigation.
\end{rem}

\section*{Acknowledgements}
We are indebted to the referee for
his/her valuable comments which have greatly improved our earlier
version.

The research of J. Bao was supported in part by the National Natural Science Foundation
of China under Grant 11401592; the research of G. Yin was supported in
part by the U.S. Army Research Office under Grant W911NF-15-1-0218; the research
of C. Yuan was supported in part by the EPSRC and NERC.


%



\printhistory
\end{document}